\documentclass[review]{elsarticle}

\usepackage{lineno,hyperref}
\modulolinenumbers[5]

\usepackage{lipsum}
\usepackage{amsfonts}
\usepackage{graphicx}
\usepackage{epstopdf}
\usepackage{algorithmic}
\ifpdf
  \DeclareGraphicsExtensions{.eps,.pdf,.png,.jpg}
\else
  \DeclareGraphicsExtensions{.eps}
\fi
\usepackage{amssymb}
\usepackage{amsthm}
\usepackage{isomath}% Pour utiliser, mathbold
\usepackage{mathtools}
\usepackage{tikz}
\usepackage{pgfplots}
\usetikzlibrary{arrows.meta}
  \pgfplotsset{compat=newest}
  \usetikzlibrary{plotmarks}
  \usepackage{grffile}

\pgfplotsset{every axis/.append style={
        scaled ticks = false, 
        tick label style={/pgf/number format/fixed}
    }
}

\usetikzlibrary{external}
\newlength\figureheight 
\newlength\figurewidth 
\usetikzlibrary{patterns}
\usetikzlibrary{plotmarks}
\usetikzlibrary{patterns}
\usepackage{tikz-3dplot}
\usepackage{multirow}
\usepackage{subfigure}
\usepackage{esint} % Pour intégrale double
\usepackage{multirow}
\usepackage{booktabs}
\usepackage{ulem}
\usepackage{cancel}
\usepackage[autostyle]{csquotes} 
\usepackage{stackengine} 
\usepackage{float}
\usepackage{setspace}
\usepackage{placeins}
\usepackage{lipsum}
\usepackage{amsfonts}
\usepackage{stmaryrd} % Pour la notation du saut
\usepackage{algorithm2e}

\newtheorem{theorem}{Theorem}
\numberwithin{theorem}{section}
\newtheorem{remark}[theorem]{Remark}

\usepackage{amsopn}

\usepackage{adjustbox}

\journal{...}

\theoremstyle{corollary}

\usepackage{lipsum}
\makeatletter
\def\ps@pprintTitle{%
 \let\@oddhead\@empty
 \let\@evenhead\@empty
 \def\@oddfoot{}%
 \let\@evenfoot\@oddfoot}
\makeatother

\begin{document}

\begin{frontmatter}

\title{High-order FDTD schemes for Maxwell's interface problems with discontinuous coefficients and complex interfaces based on the Correction Function Method}

\author{Y.-M. Law\fnref{myfootnote1}} \author{J.C. Nave\fnref{myfootnote2}}
\address{Department of Mathematics and Statistics, McGill University, Montr\'{e}al, QC, H3A 0B9, Canada.}
\fntext[myfootnote1]{yann-meing.law-kamcio@mail.mcgill.ca}
\fntext[myfootnote2]{jean-christophe.nave@mcgill.ca}

%\ead{yann-meing.law-kamcio@mcgill.ca}


\begin{abstract}
We propose high-order FDTD schemes based on the Correction Function Method (CFM) \cite{Marques2011}
	for Maxwell's interface problems with discontinuous coefficients and 
	complex interfaces. 
The key idea of the CFM is to model the correction function near an interface
	to retain the order of a finite difference approximation.
For this,
	we solve a system of PDEs based 
	on the original problem by minimizing an energy functional.
The CFM is applied to the standard Yee scheme and a fourth-order FDTD scheme. 
The proposed CFM-FDTD schemes are verified in 2-D 
	using the transverse magnetic mode (TM$_z$).
Numerical examples include scattering of magnetic and non-magnetic dielectric cylinders, 
	and problems with manufactured solutions using various complex interfaces and 
	discontinuous piecewise varying coefficients.
Long-time simulations are also performed to provide numerical evidences of the stability 
	of the proposed numerical approach.
The proposed CFM-FDTD schemes achieve up to fourth-order convergence in $L^2$-norm and 
	provide approximations devoid of spurious oscillations.
%The CFM-Yee scheme and a CFM based on a fourth-order FDTD scheme achieve respectively second and fourth 
%	order convergence in $L^2$-norm.
%Both CFM-FDTD schemes provide approximations devoid of spurious oscillations.
\end{abstract}

%\begin{keyword}
%\texttt{elsarticle.cls}\sep \LaTeX\sep Elsevier \sep template
%\MSC[2010] 00-01\sep  99-00
%\end{keyword}

\end{frontmatter}

%\linenumbers

\section{Introduction}
In computational electromagnetics, 
	the development of finite difference (FD) strategies to tackle 
	Maxwell's interface problems remains a challenge \cite{Zhang2016}.
Indeed, 
	one should expect from a numerical approach to treat arbitrary complex 
	geometries of the interface without increasing the complexity of the method,
	achieve high-order convergence to diminish the phase error for long-time simulations \cite{Hesthaven2003} and 
	handle discontinuous coefficients and discontinuous solutions,
	to name a few.
%For FD approaches, 
The potential lack of regularity of the solution of such problems is a well-known 
	challenge \cite{Leveque1994,Fedkiw1999,Marques2011}. 
%	makes difficult the use of numerical 
%	methods that require at least the solution to be continuous.
Moreover, 
	FD schemes often use simple Cartesian mesh grids and therefore
	the representation of the interface and the enforcement of interface conditions,
	fundamental to obtain accurate results, 
	are far from trivial. 
Hence,
	a first approach that consists of a staircased approximation of the interface and 
	the use of the well-known Yee scheme \cite{Yee1966}, 
	which is a second-order finite-difference time-domain (FDTD) scheme,
	yields a first-order scheme at best and non-convergent approximations 
	in some cases \cite{Ditkowski2001}.

Several numerical strategies have been proposed to overcome these issues. 
A staircase-free second-order FDTD scheme is proposed in \cite{Ditkowski2001} which explicitly enforces 
	interface conditions.
This numerical strategy has been verified for non-magnetic dielectric
	and perfect electric conductor (PEC) 
	problems using a 2-D transverse magnetic (TM) form of Maxwell's equations  \cite{Ditkowski2001,Dridi2001}. 
%In this 2-D simplification of Maxwell's equations, 
%	it is worth mentioning that electromagnetic fields are continuous across 
%	the interface between the vacuum and a non-magnetic dielectric material.
%However,
%	for a magnetic dielectric material,
%	the electric field is still continuous across the interface while the magnetic field is discontinuous.	
Inspired by the Immersed Interface Method (IIM) \cite{Leveque1994}, 
	an Upwinding Embedded Boundary (UEB) method has also been developed 
	to obtain a global second-order scheme to treat magnetic and non-magnetic dielectric problems 
	using a TM form of Maxwell's equations \cite{Cai2003}.
In the same vein, 
	high-order FDTD schemes based on the Matched Interface and Boundary (MIB) method 
	have been proposed in \cite{Zhao2004}.
These strategies derive and use jump conditions to correct a finite difference approximation in the vicinity of the interface. 
MIB-based strategies were originally limited to non-magnetic dielectrics \cite{Zhao2004,Nguyen2015} but 
	later generalized to consider a discontinuous electromagnetic field at the interface \cite{Zhang2016,Nguyen2016} using 
	2-D forms of Maxwell's equations. 
However, 
	%high-order convergence is difficult for arbitrary complex smooth interfaces. 
%Moreover, 
	the use of complex interfaces and high-order partial derivatives in jump conditions 
	increase the complexity of MIB strategies as its order increases \cite{Zhao2004,Yu2007}. 
	
Another avenue consists of FDTD schemes based on the Correction Function Method (CFM) \cite{Marques2011}.
%Like previous FDTD strategies, 
%	CFM-FDTD based strategies find a correction of a given finite difference approximation in the vicinity of 
%	the interface. 
Assuming that jumps on the interface can be smoothly extended in its vicinity, 
	the CFM models corrections that are needed to retain the order of a 
	finite difference approximation close to the interface by a system of PDEs based 
	on the original problem.
The solution of this system of PDEs is referred as the correction function.
Approximations of the correction function are then computed, 
	where it is needed, 
	by minimizing a functional which 
	is a square measure of the error associated with the correction function's system of PDEs. 
%Previous finite-difference approaches explicitly derive and impose high-order jump conditions while 
%	the CFM implicitly impose such conditions. 
%To achieve this, 
%	the CFM assumes that jumps on the interface can be smoothly extended in 
%	the vicinity of the interface 
%	and therefore requires a system of partial differential equations (PDEs) based 
%	on the original problem that models jumps in the solution.
%Afterward,
%	we define a functional that is a square measure of the error associated with the correction function's 
%	system of PDEs.
%Approximations of the jump and therefore the smooth extension of each variable are then computed, 
%	where it is needed as an approximate solution of the correction function's system of PDE, 
%	using a minimization procedure. 
Hence, 
	high-order FDTD schemes can be generated for complex interfaces 
	without significantly increasing the complexity of the proposed numerical strategy.
%However, 
%	the condition number of the matrix associated with the minimization problem increases as 
%	the order of a CFM-FDTD scheme increases \cite{LawNave2019}. 
The computational cost increases when compared with the original (i.e.\! without correction) FD scheme.
Additionally,
	a parallel implementation of the CFM can be easily performed since 
	minimization problems needed for the CFM are independent \cite{Abraham2017}.
High-order FD schemes based on the CFM have been originally developed for 
	2-D Poisson's equation with piecewise constant coefficients \cite{Marques2011,Marques2012,Marques2017} 
	as well as 3-D Poisson problems with interface jump conditions \cite{Marques2019}.
In computational electromagnetics, 
	the CFM has been extended to the wave equation \cite{Abraham2018} and 
	Maxwell's equations \cite{Marques2019} with constant coefficients.
It is also worth mentioning that high-order CFM-FDTD schemes have been proposed to handle 
	embedded PEC problems \cite{LawNave2020}.

The work presented here generalizes CFM-FDTD approaches to Maxwell's interface problems with discontinuous coefficients. 
%	in a general framework. 
We consider two FDTD schemes, 
	namely the Yee scheme and a fourth-order staggered FDTD scheme, 
	and correct them following the procedure described in \cite{LawNave2020}.
%The resulting numerical schemes are named the CFM-Yee scheme and 
%	the CFM-4$^{th}$ scheme in this work.
%In addition to scattering of non-magnetic and magnetic dielectric cylinder problems,
In addition to scattering of dielectric cylinder problems, 
	we also use problems with a manufactured solution for which complete 
	discontinuous electromagnetic 
	fields are considered to demonstrate the robustness and accuracy of the proposed numerical strategy.
%Long time simulations are performed to provide numerical evidences of the stability of the proposed 
%	CFM-FDTD schemes. % for scattering problems and problems with a manufactured solution. 
Finally, 
	we show that the correction function implicitly provides the appropriate high-order jump 
	conditions. 
Consequently, 
	high-order explicit jump conditions \cite{Zhao2004,Zhang2016} are not required for the presented method. 
	
The paper is organized as follows. 
In Section~\ref{sec:defPblm},
	we introduce Maxwell's interface problem. 
The Correction Function Method is described in Section~\ref{sec:CFM}.
In this same section, 
	we introduce the functional to be minimized based on Maxwell's equations with interface conditions. 
Then, 
	numerical examples are performed in Section~\ref{sec:numExp} to verify 
	properties of the proposed CFM-FDTD schemes. 
Finally,
	we provide conclusion and outlook in Section~\ref{sec:conclusion}.

%Simulation of problems involving Maxwell's equations with interface conditions are difficult for several reasons.
%In this work, 
%	we focus on the treatment of various complex geometries of the interface without increasing the complexity 
%	of a numerical approach for constant coefficients.
%This is the first step towards high-order finite-difference time-domain (FDTD) schemes for 
%	Maxwell's equations with interface conditions and discontinuous coefficients.	
	
\section{Definition of the Problem} \label{sec:defPblm}
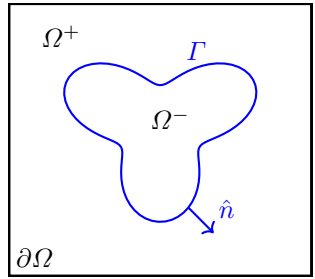
\begin{figure}[htbp]
 	\centering
		\setlength\figureheight{0.3\linewidth} 
		\setlength\figurewidth{0.35\linewidth} 
		\tikzset{external/export next=false}
\definecolor{mycolor1}{rgb}{0.00000,0.44700,0.74100}%
\begin{tikzpicture}
\begin{axis}[%
width=0.951\figurewidth,
height=\figureheight,
at={(0\figurewidth,0\figureheight)},
scale only axis,
xmin=0,
xmax=1,
ymin=0,
ymax=1,
axis line style={draw=none},
yticklabels=\empty,
xticklabels=\empty,
tick style={draw=none},
legend style={legend cell align=left,align=left,draw=white!15!black},
ylabel style={yshift=-5pt},xlabel style={yshift=2.5pt}
]
\addplot [color=blue,line width=0.8pt,solid]
  table[row sep=crcr]{%
0.75	0.55\\
0.768383691155712	0.567056285485863\\
0.784855939371776	0.586353079609766\\
0.798569247672799	0.607544505835653\\
0.808783283651588	0.630118192268234\\
0.814908287960172	0.653425855534761\\
0.816539300787274	0.676723220029148\\
0.813479080096982	0.699216694779779\\
0.805748304479489	0.720113788035931\\
0.793582461693234	0.738673991711258\\
0.777415669648107	0.754256830101566\\
0.757852505596638	0.766363942344277\\
0.735629678650627	0.774672444790369\\
0.711570021888213	0.779057373755364\\
0.686531762311539	0.779601705541283\\
0.661356319017503	0.776593244583852\\
0.636817963590651	0.770508511213973\\
0.61357854671968	0.761984593953246\\
0.592150159874138	0.75178070405503\\
0.57286808220072	0.740731832587438\\
0.555875694011952	0.729697420209933\\
0.541122262361211	0.719508274169579\\
0.528373671098689	0.710915085057808\\
0.517235331010121	0.704541799856712\\
0.507185718903124	0.700846803803962\\
0.497618308349223	0.700094370975715\\
0.48788911327206	0.702338193667119\\
0.477366703462387	0.707418035147949\\
0.465481391646248	0.714969718104235\\
0.4517703451712	0.724447815239603\\
0.435915638275791	0.735159603166056\\
0.417772716391469	0.74630812733257\\
0.397387362750012	0.757041649900452\\
0.375	0.76650635094611\\
0.351036977897187	0.773898951759353\\
0.326089339866257	0.778515940109955\\
0.300880372261704	0.779796300355625\\
0.276223968364629	0.777355071772135\\
0.252976437718672	0.771005649468394\\
0.231984821811756	0.760769465763369\\
0.214035011603472	0.746872499529074\\
0.199802985787138	0.729728904825291\\
0.189812299298019	0.709912874076284\\
0.184400561411509	0.68812060227235\\
0.183697076668585	0.66512484910402\\
0.18761311595587	0.641725065201822\\
0.195845484459602	0.61869632675674\\
0.207893209093245	0.596740392003832\\
0.223086334355692	0.576442049038377\\
0.240625045742688	0.558233576556572\\
0.259626683065719	0.542369609807538\\
0.279177704357766	0.528914027386286\\
0.298387346568558	0.517739694017157\\
0.316439622097695	0.50854106036347\\
0.332640397236881	0.500858786781107\\
0.34645661293846	0.494114777441185\\
0.357545209872456	0.487655334569043\\
0.365769976374519	0.480799613212578\\
0.371205307595373	0.472890209038588\\
0.374126697681612	0.463342567597689\\
0.374988630816851	0.451689972653066\\
0.374391337443451	0.437621149210248\\
0.37303858778221	0.421007986082703\\
0.371689260765654	0.401921513178496\\
0.37110581656491	0.380635019844539\\
0.372002990169758	0.357614024441967\\
0.375	0.33349364905389\\
0.380579330947101	0.309044762754784\\
0.389054720761967	0.285130980280279\\
0.400550380065498	0.262659193808722\\
0.414992747983782	0.242526735959632\\
0.432115274321156	0.225568494996845\\
0.45147587740097	0.212507314207483\\
0.472485908299546	0.203910805691148\\
0.494448709733373	0.200157307138778\\
0.516605239008747	0.201413134212458\\
0.538183768940384	0.207622567626085\\
0.558450417734776	0.218511208551704\\
0.576757205393503	0.233602490007808\\
0.592584493652186	0.252246299487896\\
0.605575028595216	0.273657902454885\\
0.615557346626805	0.29696470637777\\
0.622556990666661	0.321257912229455\\
0.6267947702146	0.345645796239216\\
0.628672135768096	0.369305268558684\\
0.628744571230722	0.391528473395405\\
0.627684683890353	0.411761519426597\\
0.626237340401908	0.429632939049313\\
0.625169715962851	0.444970137501007\\
0.625219459117422	0.457802865574245\\
0.627044304722356	0.468353582983459\\
0.631176384055404	0.477015419985697\\
0.637984189046328	0.484319238735193\\
0.647644665720762	0.490891992198984\\
0.660127270910301	0.497409132685517\\
0.67519106704659	0.504544198677694\\
0.692395100958555	0.512918883655448\\
0.711121467043621	0.523056852822891\\
0.730609647080229	0.535344325657581\\
0.75	0.55\\
};
\addplot [color=black,line width=1.5pt,solid]
  table[row sep=crcr]{%
1	1\\
1	0\\
0	0\\
0	1\\
1	1\\
1	1\\
};
\addplot [color=blue,line width=0.8pt,solid,->]
  table[row sep=crcr]{%
0.5926	0.2522\\
0.675925	0.159025\\
};
\end{axis}
\draw (2.5,3.0) node {\color{blue}$\Gamma$};	
\draw (2.9,0.9) node {\color{blue}$\hat{n}$};
\draw (2.15,2.1) node {\color{black}$\Omega^-$};
\draw (0.7,3.2) node {\color{black}$\Omega^+$};
\draw (0.35,0.25) node {\color{black}$\partial\Omega$};
\end{tikzpicture}%  
		\caption{Geometry of a domain $\Omega$ with an interface $\Gamma$.}
\label{fig:typicalDomain}
\end{figure}

Assume a domain in space $\Omega$ subdivided 
	into two subdomains $\Omega^+$ and $\Omega^-$,
 	and a time interval $I = [0,T]$.
 The interface $\Gamma$ between subdomains is independent of time and allows the solutions to 
 	be discontinuous. 
Figure~\ref{fig:typicalDomain} illustrates a typical geometry of a domain $\Omega$.
%%%%%%%%%%%%%%%%%%%%%%%%%%%%%%%%%%%%%%%%%%
%%%%%%%%%%%%%%%%%%%%%%%%%%%%%%%%%%%%%%%%%%	
For a given variable $\mathbold{A}$, 
	we define $\mathbold{A}^+$ and $\mathbold{A}^-$ as respectively the solutions in $\Omega^+$ and $\Omega^-$.
A jump of $\mathbold{A}$ on the interface $\Gamma$ is denoted as $\llbracket \mathbold{A} \rrbracket = \mathbold{A}^+ - \mathbold{A}^-$. 
Assuming linear media, 
	we consider Maxwell's equations with interface conditions that are given by 
\normalsize
\begin{subequations} \label{eq:pblmDefinition}
\begin{align}
\mu(\mathbold{x})\,\partial_t \mathbold{H} + \nabla\times \mathbold{E} =&\,\, 0 \quad \text{in } \Omega \times I, \label{eq:Faraday} \\
\epsilon(\mathbold{x})\,\partial_t \mathbold{E} - \nabla\times\mathbold{H} =&\,\, 0 \quad \text{in } \Omega \times I , \label{eq:AmpereMaxwell}\\
\nabla\cdot(\epsilon(\mathbold{x})\,\mathbold{E}) =&\,\, 0 \quad \text{in } \Omega \times I , \label{eq:divD}\\
\nabla\cdot(\mu(\mathbold{x})\,\mathbold{H})=&\,\, 0 \quad \text{in } \Omega \times I , \label{eq:divB}\\
\hat{\mathbold{n}}\times\llbracket \mathbold{E} \rrbracket =&\,\, 0 \quad \text{on } \Gamma \times I ,\label{eq:tangentEInterf}\\
\hat{\mathbold{n}}\times\llbracket \mathbold{H} \rrbracket =&\,\, 0  \quad \text{on } \Gamma \times I ,\label{eq:tangentHInterf}\\
\hat{\mathbold{n}}\cdot\llbracket \epsilon(\mathbold{x})\,\mathbold{E} \rrbracket =&\,\, 0 \quad \text{on } \Gamma \times I ,\label{eq:normalEInterf}\\
\hat{\mathbold{n}}\cdot\llbracket \mu(\mathbold{x})\,\mathbold{H} \rrbracket =&\,\, 0 \quad \text{on } \Gamma \times I ,\label{eq:normalHInterf}\\
\mathbold{n}\times\mathbold{H} =&\,\, \mathbold{g}_1(\mathbold{x},t)	\quad \text{on } \partial\Omega\times I, \label{eq:bndCdnH}\\
\mathbold{n}\times\mathbold{E} =&\,\, \mathbold{g}_2(\mathbold{x},t)	\quad \text{on } \partial\Omega \times I, \label{eq:bndCdnE}\\
\mathbold{H} =&\,\, \mathbold{H}(\mathbold{x},0)	\quad \text{in } \Omega, \label{eq:InitialCdnH}\\
\mathbold{E} =&\,\, \mathbold{E}(\mathbold{x},0)	\quad \text{in } \Omega, \label{eq:InitialCdnE}
\end{align}
\end{subequations}
\normalsize
	where $\mathbold{H}$ is the magnetic field, 
	$\mathbold{E}$ is the electric field, 
	$\mu(\mathbold{x})>0$ is the magnetic permeability, 
	$\epsilon(\mathbold{x})>0$ is the electrical permittivity,
%	$\rho$ is the electric charge density,
%	$\mathbold{J}_s$ is the surface current density, 
%	$\rho_s$ is the surface charge density,
	$\mathbold{n}$ is the unit outward normal to $\partial\Omega$  and 
	$\hat{\mathbold{n}}$ is the unit normal to the interface $\Gamma$ pointing toward $\Omega^+$.
Interface conditions are given by equations \eqref{eq:tangentEInterf} to \eqref{eq:normalHInterf} 
	while boundary and initial conditions are given by equations \eqref{eq:bndCdnH} 
	to \eqref{eq:InitialCdnE}.
Physical parameters, 
	that is $\mu$ and $\epsilon$, 
	can be discontinuous on the interface.
%In this work, 
Without loss of generality, 
	we assume that electromagnetic fields 
%	that is $\mathbold{H}$ and $\mathbold{E}$,
	are at divergence-free in each subdomain.

\section{Correction Function Method} \label{sec:CFM}
The Correction Function Method (CFM) allows one to find a correction for a given FD approximation 
	involving nodes that belong to different subdomains in order to retain its order.
 For this purpose,
 	the CFM assumes that solutions in each subdomain can be extended across 
	the interface $\Gamma$ in a small domain $\Omega_\Gamma\times I$, 
	that is such that $\Omega_\Gamma \subset \Omega$ encloses $\Gamma$.
%Afterward, 
%	a correction of a finite difference approximation that involves nodes in different subdomains,
%	that is on both sides of the interface $\Gamma$, 
%	can be found to retain the order of the FD scheme. 
A system of PDEs based on the original problem, 
	namely Maxwell's interface problem \eqref{eq:pblmDefinition} in our case,  
	models the extension of each variable around the interface. % while satisfying the 
%	interface conditions. 
The solution of this system of PDEs is referred as the correction function. 
Afterward, 
	we define a functional that is a square measure of the error associated with the correction function's 
	system of PDEs.
Approximations of the correction function are then computed, 
	where it is needed, 
	using a minimization procedure. 
In practice, 
	the interface is discretized and a local patch
%	$\Omega_\Gamma^h\subset \Omega_\Gamma$
	$\Omega_\Gamma^h\times I_\Gamma^h \subset \Omega_\Gamma \times I$
	is defined for each node of its discretization. 
Moreover, 
	the size of local patches depends on the considered FD scheme and % original FD scheme, that is without correction,
	should diminish as the mesh grid size diminishes (see Remark~\ref{rem:sizePatch}).
%%The time interval of local patches $I_\Gamma^h$ depends on the considered time-stepping method.
%%As for the size of the domain in space $\Omega_\Gamma^h$, 
%%	it is chosen in such a way that each node $(\mathbold{x}_c,t_n)$ to be corrected should be within a local patch 
%%	and its size should diminish as the mesh grid size diminishes.
%We refer to \cite{LawNave2020} for more details about the construction of local patches and 
%	the application of the CFM on FDTD staggered schemes.
%Afterward, 
%	each node $(\mathbold{x}_c,t_n)$ to be corrected is assigned to a local patch.
%It is worth mentioning that the additional cost associated with this approach is not negligible and 
%	consumes most of the computational time when compared with the original FD scheme, 
%	that is without correction. 
%However,  
%	a parallel implementation of the CFM greatly reduces its computational time and
%	 can be easily performed since minimization problems for a given time-step are independent \cite{Abraham2017}. 

In the following, 
	we derive the system of PDEs that models the smooth extension of each variable 
	and therefore the correction function. % using Maxwell's equations with interface conditions \eqref{eq:pblmDefinition}.
The minimization problem based on the associated energy functional is also presented.
%	a functional that is a square measure of the error associated with the 
%	correction function's system of PDEs is presented. 
%Afterward, 
%	we analyze the minimization problem and identify under which conditions it has a 
%	global minimizer. 

Let us first introduce some notations.	
The inner product in $L^2\big(\Omega_{\Gamma}^{h}\times I_{\Gamma}^h\big)$ is defined by
$$\langle\mathbold{v},\mathbold{w}\rangle = \int\limits_{ I_{\Gamma}^h}\!\int\limits_{\Omega_{\Gamma}^{h}}\!\!\mathbold{v}\cdot\mathbold{w}\,\mathrm{d}V\,\mathrm{d}t$$
with $\|\mathbold{v}\| = \langle\mathbold{v},\mathbold{v}\rangle$,
and we also use the notation
	$$\langle \mathbold{v},\mathbold{w} \rangle_{\Gamma} = \int\limits_{ I_{\Gamma}^h}\!\int\limits_{\Gamma\cap\Omega_{\Gamma}^{h}}\!\!\mathbold{v}\cdot\mathbold{w} \, \mathrm{d}S\,\mathrm{d}t$$
	with $\|\mathbold{v}\|_{\Gamma} = \langle\mathbold{v},\mathbold{v}\rangle_{\Gamma}$ for legibility. 
%Unlike previous CFM-FDTD schemes for Maxwell's equations, 
Unlike previous CFM-FDTD schemes, 
%Unlike previous CFM-FD schemes, 
	we cannot explicitly model jumps $\mathbold{D}_H = \llbracket \mathbold{H} \rrbracket$ 
	and $\mathbold{D}_E = \llbracket \mathbold{E} \rrbracket$ because of discontinuous coefficients.
Hence, 
	we first need to estimate $\mathbold{H}^+$, 
	$\mathbold{E}^+$,
	$\mathbold{H}^-$ and $\mathbold{E}^-$ in the whole patch, 
	and afterward compute an approximation of $\mathbold{D}_H$ and $\mathbold{D}_E$. 
The system of PDEs for correction functions is then given by  
\begin{equation} \label{eq:pblmDefinitionCFM}
\begin{aligned}
\mu^+(\mathbold{x})\,\partial_t \mathbold{H}^+ + \nabla\times \mathbold{E}^+ =&\,\, 0 \quad \text{in } \Omega_\Gamma^h\times I_\Gamma^h ,  \\
\epsilon^+(\mathbold{x})\,\partial_t \mathbold{E}^+ - \nabla\times\mathbold{H}^+ =&\,\, 0 \quad \text{in } \Omega_\Gamma^h\times I_\Gamma^h ,\\
\nabla\cdot(\epsilon^+(\mathbold{x})\,\mathbold{E}^+) =&\,\, 0 \quad \text{in } \Omega_\Gamma^h\times I_\Gamma^h , \\
\nabla\cdot(\mu^+(\mathbold{x})\,\mathbold{H}^+)=&\,\, 0 \quad \text{in } \Omega_\Gamma^h\times I_\Gamma^h , \\
\mu^-(\mathbold{x})\,\partial_t \mathbold{H}^- + \nabla\times \mathbold{E}^- =&\,\, 0 \quad \text{in } \Omega_\Gamma^h\times I_\Gamma^h,  \\
\epsilon^-(\mathbold{x})\,\partial_t \mathbold{E}^- - \nabla\times\mathbold{H}^- =&\,\, 0 \quad \text{in } \Omega_\Gamma^h\times I_\Gamma^h ,\\
\nabla\cdot(\epsilon^-(\mathbold{x})\,\mathbold{E}^-) =&\,\, 0 \quad \text{in } \Omega_\Gamma^h\times I_\Gamma^h , \\
\nabla\cdot(\mu^-(\mathbold{x})\,\mathbold{H}^-)=&\,\, 0 \quad \text{in } \Omega_\Gamma^h\times I_\Gamma^h , \\
\hat{\mathbold{n}}\times\llbracket \mathbold{E} \rrbracket =&\,\, 0 \quad \text{on } \Gamma\cap\Omega_{\Gamma}^{h} \times I_\Gamma^h ,\\
\hat{\mathbold{n}}\times\llbracket \mathbold{H} \rrbracket =&\,\, 0  \quad \text{on } \Gamma\cap\Omega_{\Gamma}^{h} \times I_\Gamma^h ,\\
\hat{\mathbold{n}}\cdot\llbracket \epsilon(\mathbold{x})\,\mathbold{E} \rrbracket =&\,\, 0 \quad \text{on } \Gamma\cap\Omega_{\Gamma}^{h} \times I_\Gamma^h ,\\
\hat{\mathbold{n}}\cdot\llbracket \mu(\mathbold{x})\,\mathbold{H} \rrbracket =&\,\, 0 \quad \text{on } \Gamma\cap\Omega_{\Gamma}^{h} \times I_\Gamma^h ,
\end{aligned}
\end{equation}
Following the procedure described in \cite{LawMarquesNave2020} to construct 
	a functional that is a square measure of the error associated with 
	system \eqref{eq:pblmDefinitionCFM}
	leads to an ill-posed minimization problem.  
%	, one obtains an ill-posed minimization problem. % in the limit when the mesh grid size goes to zero.
% This could lead to ill-conditioned matrices coming from the minimization problem.
As in CFM-FDTD strategies for embedded perfect electric conductors \cite{LawNave2020}, 
	we can take advantage of FD approximations at previous time steps using fictitious interface 
	conditions to retrieve a well-posed minimization problem.
%     reduce the condition number of these matrices. 
%	retrieve 
%	a well-posed minimization problem (c.f. Corollary~\ref{cor:minGlobal}). 
Fictitious interface conditions are given by 
	\begin{equation} \label{eq:fictInterfCdn}
	\begin{aligned}
	\hat{\mathbold{n}}_{1,i}^\circ\times(\mathbold{E}^\circ-\mathbold{E}^{\circ,*}) =&\,\, 0 \quad \text{on} \quad \Gamma_{1,i}^\circ \times I_\Gamma^h  \quad \text{for} \quad i=1,\ldots,N_1^\circ,\\
	\hat{\mathbold{n}}_{2,i}^\circ\times(\mathbold{H}^\circ-\mathbold{H}^{\circ,*}) =&\,\, 0  \quad \text{on} \quad \Gamma_{2,i}^\circ \times I_\Gamma^h \quad \text{for} \quad i=1,\ldots,N_2^\circ,\\
	\hat{\mathbold{n}}_{3,i}^\circ\cdot(\mathbold{E}^\circ-\mathbold{E}^{\circ,*}) =&\,\, 0 \quad \text{on} \quad \Gamma_{3,i}^\circ \times I_\Gamma^h \quad \text{for} \quad i=1,\ldots,N_3^\circ,\\
	\hat{\mathbold{n}}_{4,i}^\circ\cdot(\mathbold{H}^\circ-\mathbold{H}^{\circ,*}) =&\,\, 0 \quad \text{on} \quad \Gamma_{4,i}^\circ \times I_\Gamma^h \quad \text{for} \quad i=1,\ldots,N_4^\circ,
\end{aligned}
\end{equation}
	where $\circ$ is either $+$ or $-$ depending in which subdomain the fictitious interface 
	$\Gamma_{k,i}^\circ$ belongs, 
	$\hat{\mathbold{n}}_{k,i}^\circ$ is the normal associated with $\Gamma_{k,i}^\circ$,
	$N_k^\circ$ is the number of fictitious interfaces,
	and $\mathbold{H}^{\circ,*}$ and $\mathbold{E}^{\circ,*}$ are approximations of the magnetic field 
	and the electric field that come from a FD scheme.

The quadratic functional to minimize is therefore given by
\small
\begin{equation*}
\begin{aligned}
J(\mathbold{H}^+,&\mathbold{E}^+,\mathbold{H}^-,\mathbold{E}^-) = \frac{\ell_h}{2} \, \big\|\mu^+\,\partial_t \mathbold{H}^+ + \nabla\times\mathbold{E}^+\big\| +  \frac{\ell_h}{2} \, \big\| \epsilon^+\,\partial_t \mathbold{E}^+ - \nabla\times\mathbold{H}^+\big\| \\
+&\,\,  \frac{\ell_h}{2} \, \big\| \nabla\cdot(\epsilon^+\, \mathbold{E}^+)\big\| +  \frac{\ell_h}{2} \, \big\| \nabla\cdot( \mu^+\, \mathbold{H}^+ )\big\|  + \frac{\ell_h}{2} \, \big\|\mu^-\,\partial_t \mathbold{H}^- + \nabla\times\mathbold{E}^-\big\| \\
+&\,\,  \frac{\ell_h}{2} \, \big\| \epsilon^-\,\partial_t \mathbold{E}^- - \nabla\times\mathbold{H}^-\big\| + \frac{\ell_h}{2} \, \big\| \nabla\cdot(\epsilon^-\, \mathbold{E}^-)\big\| +  \frac{\ell_h}{2} \, \big\| \nabla\cdot( \mu^-\, \mathbold{H}^- )\big\| \\
+&\,\, \frac{c_p}{2} \, \big\|\hat{\mathbold{n}}\times(\mathbold{E}^+-\mathbold{E}^-)\big\|_{\Gamma} + \frac{c_p}{2} \,\big\| \hat{\mathbold{n}}\times(\mathbold{H}^+-\mathbold{H}^-)\big\|_{\Gamma} \\
+&\,\, \frac{c_p}{2} \, \big\|\hat{\mathbold{n}}\cdot(\epsilon^+\,\mathbold{E}^+-\epsilon^-\,\mathbold{E}^-) \big\|_{\Gamma} + \frac{c_p}{2}\,\big\|\hat{\mathbold{n}}\cdot(\mu^+\,\mathbold{H}^+-\mu^-\,\mathbold{H}^-)\big\|_{\Gamma}\\
+&\,\, \frac{c_f}{2\,N_{E^+}} \, \sum_{i=1}^{N_1^+} \big\| \hat{\mathbold{n}}_{1,i}^+\times(\mathbold{E}^+-\mathbold{E}^{+,*})\big\|_{\Gamma_{1,i}^+} 
+ \frac{c_f}{2\,N_{H^+}} \, \sum_{i=1}^{N_2^+} \big\| \hat{\mathbold{n}}_{2,i}^+\times(\mathbold{H}^+-\mathbold{H}^{+,*})\big\|_{\Gamma_{2,i}^+} \\
+&\,\, \frac{c_f}{2\,N_{E^+}} \, \sum_{i=1}^{N_3^+}  \big\|\hat{\mathbold{n}}_{3,i}^+\cdot(\mathbold{E}^+-\mathbold{E}^{+,*}) \big\|_{\Gamma_{3,i}^+} 
+ \frac{c_f}{2\,N_{H^+}} \, \sum_{i=1}^{N_4^+}  \big\|\hat{\mathbold{n}}_{4,i}^+\cdot (\mathbold{H}^+-\mathbold{H}^{+,*}) \big\|_{\Gamma_{4,i}^+}\\
+&\,\, \frac{c_f}{2\,N_{E^-}} \, \sum_{i=1}^{N_1^-} \big\| \hat{\mathbold{n}}_{1,i}^-\times(\mathbold{E}^--\mathbold{E}^{-,*})\big\|_{\Gamma_{1,i}^-} 
+ \frac{c_f}{2\,N_{H^-}} \, \sum_{i=1}^{N_2^-} \big\| \hat{\mathbold{n}}_{2,i}^-\times(\mathbold{H}^--\mathbold{H}^{-,*})\big\|_{\Gamma_{2,i}^-} \\
+&\,\, \frac{c_f}{2\,N_{E^-}} \, \sum_{i=1}^{N_3^-}  \big\|\hat{\mathbold{n}}_{3,i}^-\cdot(\mathbold{E}^--\mathbold{E}^{-,*}) \big\|_{\Gamma_{3,i}^-} 
+ \frac{c_f}{2\,N_{H^-}} \, \sum_{i=1}^{N_4^-}  \big\|\hat{\mathbold{n}}_{4,i}^-\cdot (\mathbold{H}^--\mathbold{H}^{-,*}) \big\|_{\Gamma_{4,i}^-}\\
\end{aligned}
\end{equation*}	
\normalsize
	where $c_p>0$ and $c_f>0$ are penalization coefficient, 
	$\ell_h$ is the characteristic length in space of the patch,
	$N_{E^\circ} = N_1^\circ + N_3^\circ$ and 
	$N_{H^\circ} = N_2^\circ + N_4^\circ$.
Integrals over the domain are scaled by $\ell_h$ to guarantee that all terms 
	in the functional $J$ behave in a similar way when the 
	computational grid is refined \cite{LawMarquesNave2020}.
The problem statement is then 
\begin{equation} \label{eq:minPblm}
\begin{aligned}
&\text{Find } (\mathbold{H}^+,\mathbold{E}^+,\mathbold{H}^-,\mathbold{E}^-) \in V \times W \times V \times W \text{ such that }\\
 &\qquad (\mathbold{H}^+,\mathbold{E}^+,\mathbold{H}^-,\mathbold{E}^-) \in  \underset{\substack{\mathbold{v}^+,\mathbold{v}^- \in V \\ \mathbold{w}^+,\mathbold{w}^- \in W}}{\arg\min}J(\mathbold{v}^+,\mathbold{w}^+,\mathbold{v}^-,\mathbold{w}^-),
\end{aligned}
\end{equation}
	where $W=V$.
Let us recall that we assume divergence-free electromagnetic fields in each subdomain.
We therefore minimize the functional $J$ in a space of divergence-free space-time polynomials,
	namely
\begin{equation*}
	V = \big\{ \mathbold{v} \in \big[P^k\big(\Omega_{\Gamma}^{h}\times I^h_{\Gamma}\big)\big]^3 : \nabla\cdot\mathbold{v} = 0 \big\},
\end{equation*}
	where $P^k$ denotes the space of polynomials of degree $k$.
It is worth mentioning that basis functions of $V$ are based on high-degree divergence-free 
	basis functions proposed in \cite{Cockburn2004}.
%The degree $k$ of polynomial basis functions is chosen to be at least equal to the order of 
%	a chosen FDTD scheme to retain its order.
%It is also worth to mention that we also use space-time polynomial interpolations for FD approximations 
%	of at least the same order of the considered CFM-FDTD scheme to compute fictitious 
%	interface conditions.
\begin{remark} \label{rem:sizePatch}
The size in space of local patches $\ell_h$ depends of the mesh grid size, 
	that is $\ell_h = \beta\,\max\{\Delta x, \Delta y, \Delta z\}$,
	where $\beta$ is a positive constant. 
The choice of $\beta$ depends on the considered FD scheme and must allow 
	the construction of enough fictitious interfaces within the local patch.
To ease the implementation, 
	local patches are taken aligned with the mesh grid and square in space. 
Fictitious interfaces are also aligned with the mesh grid to facilitate the computation of space-time interpolants that are 
	needed in the minimization problem. 
%Hence,
% 	we can directly use FD approximations to compute space-time interpolants that are 
%	needed in the minimization problem. 
We refer the reader to \cite{LawNave2020} for more details on the implementation 
	of local patches and fictitious interface conditions.
\end{remark}

\begin{remark} \label{rem:initializationScheme}
The initialization of CFM-FDTD schemes can be difficult because of time integrals involving 
	$\mathbold{H}^*$ and $\mathbold{E}^*$. 
An initialization strategy has been developed for the 
	Yee scheme and a fourth-order FDTD scheme based on a 
	multistep method \cite{LawNave2020}.
Another approach, 
	which is specific to some applications,
%	such as radar,
	consists to assume that electromagnetic fields close to the interface 
	remain unchanged for $t\leq t_0$.	
\end{remark}

\begin{remark}
Using a truncation error analysis, 
	one can show that the order of a CFM-FDTD scheme for Maxwell's 
	equations \eqref{eq:pblmDefinition} is $\min\{n,k\}$ where $n$ is the order of
	the considered FD scheme and $k$ is the degree of space-time polynomial spaces used in minimization problem \eqref{eq:minPblm} \cite{LawNave2020}.
\end{remark}
\begin{remark} \label{rem:jumpCdns}
%The correction function's system of PDEs on which the functional $J$ is based models 
%	the behaviour of jumps in electromagnetic fields in the vicinity of the interface.
The correction function's system of PDEs on which functional $J$ is based models 
	the extension of each electromagnetic field in the vicinity of the interface while 
	satisfying interface conditions.
Hence, 
	by construction and consistency, 
	explicit jump conditions on the interface used for Matched Interface and Boundary based strategies 
	\cite{Zhao2004,Zhang2016} should be implicitly satisfied.
This claim is supported by numerical evidences presented in 
	subsection~\ref{sec:numExpScattering}.
\end{remark}

\begin{remark} \label{rem:stability}
It is recalled that fictitious interface conditions are used to retrieve a well-posed minimization problem.
%	and reduce the condition number 
%	of matrices coming from the minimization procedure. 
Regarding the value of $c_f$, 
	the priority should be given to interface conditions and therefore $c_p>c_f>0$.
Moreover, 
	$c_f$ should also diminish as the mesh grid size diminishes to enforce again interface conditions.
As mentioned in \cite{LawNave2020}, 
	the stability analysis of a CFM-FDTD scheme that uses 
	fictitious interface conditions \eqref{eq:fictInterfCdn} is not trivial.
Despite the lack of a rigorous proof, 
	$c_f = \alpha \, \Delta t$,
	where $\alpha$ is a positive constant sufficiently small,
	seems to avoid any stability issues.
We also assume that the stability condition of a CFM-FDTD scheme should be close to the one 
	associated with the original (i.e.\! without correction) FDTD scheme.
This is corroborated with numerical results in \cite{LawNave2020} and the 
	performed numerical examples in subsection~\ref{sec:longTimeSimulations}.
\end{remark}

%\subsection{Fictitious interface conditions for the CFM} \label{sec:fictInterfCFM}
%	
%\subsection{Computation of the local patch} \label{sec:compLocalPatch}
%
%\subsubsection{Stability Analysis}

%%%%%%%%%%%%%%%%%%%%%%%%%%%%%%%%%%%%%%%%%%%%%%%%%
\FloatBarrier
\section{Numerical Examples} \label{sec:numExp}
In this section, 
	we perform convergence analysis and long-time simulations in 2-D to verify 
	the proposed numerical strategy. 
We consider the transverse magnetic (TM$_{\text{z}}$) mode.
Hence,
	for a domain $\Omega \subset \mathbb{R}^2$, 
	Maxwell's equations are simplified to
\begin{equation} \label{eq:TMzSyst}
	\begin{aligned}
	\mu(x,y)\,\partial_t H_x + \partial_y E_z =&\,\, 0 \quad \text{in } \Omega \times I,\\
	\mu(x,y)\,\partial_t H_y - \partial_x E_z  =&\,\, 0 \quad \text{in } \Omega \times I,\\
	\epsilon(x,y)\,\partial_t E_z - \partial_x H_y + \partial_y H_x =&\,\, 0 \quad \text{in } \Omega \times I,\\
	\partial_x (\mu(x,y)\,H_x) + \partial_y (\mu(x,y)\,H_y) =&\,\, 0 \quad \text{in } \Omega \times I,\\
	\end{aligned}
\end{equation}
	with the associated interface, 
	boundary and initial conditions. 
In this 2-D simplification of Maxwell's equations, 
	electromagnetic fields are continuous across 
	the interface between the vacuum and a non-magnetic dielectric material.
However,
	for a magnetic dielectric material,
	the electric field is still continuous across the interface while the magnetic field is discontinuous.
	
We consider two different FDTD schemes,
	namely the Yee scheme and a fourth-order FDTD scheme.
The latter FDTD scheme also uses staggered grids in both space and time.
More specifically,
	space derivatives are estimated with a fourth-order centered 
	FD approximation for staggered grids while time derivatives are estimated using 
	a fourth-order staggered free-parameter multistep method \cite{Ghrist2000}.
The associated CFM-FDTD schemes are then the CFM-Yee scheme and the CFM-4$^{th}$ scheme.
We refer to \cite{LawNave2020} for more details on the application of the CFM to these two FDTD 
	schemes.

\subsection{Scattering of a Dielectric Cylinder Problems} \label{sec:numExpScattering}
Let us consider a dielectric cylinder in free-space exposed to a TM$_z$ excitation wave. 
The interface is a circle of radius $r_0=0.6$ centered at $(0,0)$. 
The exact solution in cylindrical coordinates is given by the real part of 
\footnotesize
\begin{equation*}
\begin{aligned}
H_{\theta}(r,\theta,t) =&\,\, \left\{ 
  \begin{array}{l l}
    -\frac{\mathfrak{i}\,k^-}{\omega\,\mu^-}\displaystyle\sum_{n=-\infty}^{\infty} C_n^{\text{tot}}\,J_n^\prime(k^-\,r)\,e^{\mathfrak{i}\,(n\,\theta+\omega\,t)}, & \text{if} \quad r\leq r_0, \\ %[1.5mm]
    -\frac{\mathfrak{i}\,k^+}{\omega\,\mu^+}\displaystyle\sum_{n=-\infty}^{\infty} (\mathfrak{i}^{-n}\,J_n^\prime(k^+\,r) + C_n^{\text{scat}}\,H_n^{{(2)}^\prime}(k^+\,r))\,e^{\mathfrak{i}\,(n\,\theta+\omega\,t)}, & \text{if} \quad r>r_0,
  \end{array} \right.\\
H_r(r,\theta,t) =&\,\, \left\{ 
  \begin{array}{l l}
    -\frac{1}{\omega\,\mu^-\,r}\displaystyle\sum_{n=-\infty}^{\infty} n\,C_n^{\text{tot}}\,J_n(k^-\,r)\,e^{\mathfrak{i}\,(n\,\theta+\omega\,t)}, & \text{if} \quad r\leq r_0, \\ %[1.5mm]
    -\frac{1}{\omega\,\mu^+\,r}\displaystyle\sum_{n=-\infty}^{\infty} n\,(\mathfrak{i}^{-n}\,J_n(k^+\,r) + C_n^{\text{scat}}\,H_n^{(2)}(k^+\,r))\,e^{\mathfrak{i}\,(n\,\theta+\omega\,t)}, & \text{if} \quad r>r_0,
  \end{array} \right. \\
E_z(r,\theta,t) =&\,\, \left\{ 
  \begin{array}{l l}
    \displaystyle\sum_{n=-\infty}^{\infty} C_n^{\text{tot}}\,J_n(k^-\,r)\,e^{\mathfrak{i}\,(n\,\theta+\omega\,t)}, & \text{if} \quad r\leq r_0, \\ %[1.5mm]
    \displaystyle\sum_{n=-\infty}^{\infty} (\mathfrak{i}^{-n}\,J_n(k^+\,r) + C_n^{\text{scat}}\,H_n^{(2)}(k^+\,r))\,e^{\mathfrak{i}\,(n\,\theta+\omega\,t)}, &  \text{if} \quad r>r_0,
  \end{array} \right.
  \end{aligned}
\end{equation*}
\normalsize
with 
\begin{equation*}
\begin{aligned}
C_n^{\text{tot}} =&\,\, \mathfrak{i}^{-n}\,\frac{\tfrac{k^+}{\mu^+}\,(J_n^\prime(k^+\,r_0)\,H_n^{(2)}(k^+\,r_0)-H_n^{{(2)}^\prime}(k^+\,r_0)\,J_n(k^+\,r_0))}{\frac{k^-}{\mu^-}\,J_n^\prime(k^-\,r_0)\,H_n^{(2)}(k^+\,r_0)-\tfrac{k^+}{\mu^+}\,H_n^{{(2)}^\prime}(k^+\,r_0)\,J_n(k^-\,r_0)},\\
C_n^{\text{scat}}=&\,\, \mathfrak{i}^{-n}\,\frac{\tfrac{k^+}{\mu^+}\,J_n^\prime(k^+\,r_0)\,J_n(k^-\,r_0)-\tfrac{k^-}{\mu^-}\,J_n^\prime(k^-\,r_0)\,J_n(k^+\,r_0)}{\frac{k^-}{\mu^-}\,J_n^\prime(k^-\,r_0)\,H_n^{(2)}(k^+\,r_0)-\tfrac{k^+}{\mu^+}\,H_n^{{(2)}^\prime}(k^+\,r_0)\,J_n(k^-\,r_0)},
\end{aligned}
\end{equation*}
where $\mathfrak{i}$ is the imaginary number, 
	$k^\circ=\omega\,\sqrt{\mu^\circ\,\epsilon^\circ}$, 
	$\omega = 2\,\pi$,
	$J_n$ is the $n$-order Bessel function of first kind and $H_n^{(2)}$ is the $n$-order 
	Hankel function of second kind \cite{Taflove1993,Cai2003}. 

\begin{figure}[htbp]
 	\centering
 	\tdplotsetmaincoords{75}{105}
	\tikzset{external/export next=false}
  \begin{tikzpicture}[scale=0.75]
   	\draw[-latex,thick] (-0.5,-0.5)--(-0.5,5.5)--(5.5,5.5)--(5.5,-0.5)--cycle; 
   	\draw[-latex,thick,dashed,magenta] (0,0)--(0,5)--(5,5)--(5,0)--cycle; 
	\draw[-latex,thick,blue] (2.5,2.5) circle [radius=1.5];
  	\draw (3.8,3.9) node {$\color{blue}\Gamma$};
  	\draw (0.8,4.2) node {$\Omega^+$};
   	\draw (2.5,2.5) node {$\Omega^-$};
  	\draw (0.5,0.3) node {\color{magenta}$\partial \Omega$};	
  	\draw (-1.1,-0.5) node {$\partial \Omega_c$};
%	\draw[-latex,thick,blue] (3.56066,1.43934)--(4.26,0.73);
%  	\draw (4.1,1.4) node {$\color{blue}\hat{\mathbold{n}}$};
%	\draw[-latex,thick,blue] (2.5,0)--(2.5,-1);
%  	\draw (4.1,1.4) node {$\color{magenta}{\mathbold{n}}$};
  \end{tikzpicture}
  \caption{Computational domain of scattering of a dielectric cylinder problems.}
\label{fig:computationalDomain}
\end{figure}
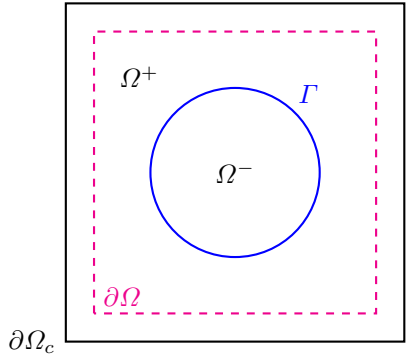

For the CFM-Yee scheme, 
	the domain is $\Omega = [-1,1]\times[-1,1]$ and we impose Dirichlet boundary conditions on 
	the boundary  $\partial \Omega$ of the domain. 
As for the CFM-4$^{th}$ scheme,  
	the domain $\Omega = [-0.9,0.9]\times[-0.9,0.9]$ is embedded in a computational domain,
	namely $\Omega_c = [-1,1]\times[-1,1]$, 
	as illustrated in Figure~\ref{fig:computationalDomain}.
%%%%%%%%%%%%%%%%%%%%%%%%%%%%%%%%%%%%%%%%%%
%%%%%%%%%%%%%%%%%%%%%%%%%%%%%%%%%%%%%%%%%%	
We use the CFM with constant coefficients to enforce electromagnetic fields on $\partial \Omega$ \cite{LawMarquesNave2020}.
%	, that is $\mathbold{n}\times\mathbold{H}$, 
%	$\mathbold{n}\cdot\mathbold{H}$,
%	$\mathbold{n}\times\mathbold{E}$ and $\mathbold{n}\cdot\mathbold{E}$,
%	where $\mathbold{n}$ is the unit outward normal 
%	to $\partial \Omega$ \cite{LawMarquesNave2020}.
Hence,
 	the trivial solution is imposed in $\Omega_c\backslash\Omega$ and 
	periodic conditions are imposed on $\partial \Omega_c$. 
The time interval is $I=[0,1]$.
The mesh grid size is $h = \Delta x = \Delta y$ with 
	$h \in \big\{ \tfrac{1}{20}, \tfrac{1}{28}, \tfrac{1}{40},\tfrac{1}{52}, \tfrac{1}{72}, \tfrac{1}{96}, \tfrac{1}{132}, \tfrac{1}{180},\tfrac{1}{244},\tfrac{1}{336}\big\}$ and the time step is $\Delta t = \tfrac{h}{2}$.
For both schemes, 
	we choose $\ell_h = 7\,h$ to construct local patches and 
	we use at least a second degree interpolating polynomial in space to construct 
	$\mathbold{H}^*$ and $\mathbold{E}^*$ that are needed for 
	fictitious interface conditions \eqref{eq:fictInterfCdn}.	
We set $c_f =\Delta t$ and $c_f = \tfrac{\Delta t}{4}$ for respectively the CFM-Yee and 
	the CFM-$4^{th}$ scheme while $c_p = 1$ for both schemes.
Second and third degree polynomial approximations of correction functions are chosen for 
	respectively the CFM-Yee and the CFM-$4^{th}$ scheme.

\begin{figure}     
	\centering
	\subfigure[a non-magnetic dielectric material ($\mu^+=\mu^-=1$, $\epsilon^+=1$ and $\epsilon^-=2.25$)]{
	\stackunder[5pt]{
	\stackunder[5pt]{\includegraphics[width=1.58in]{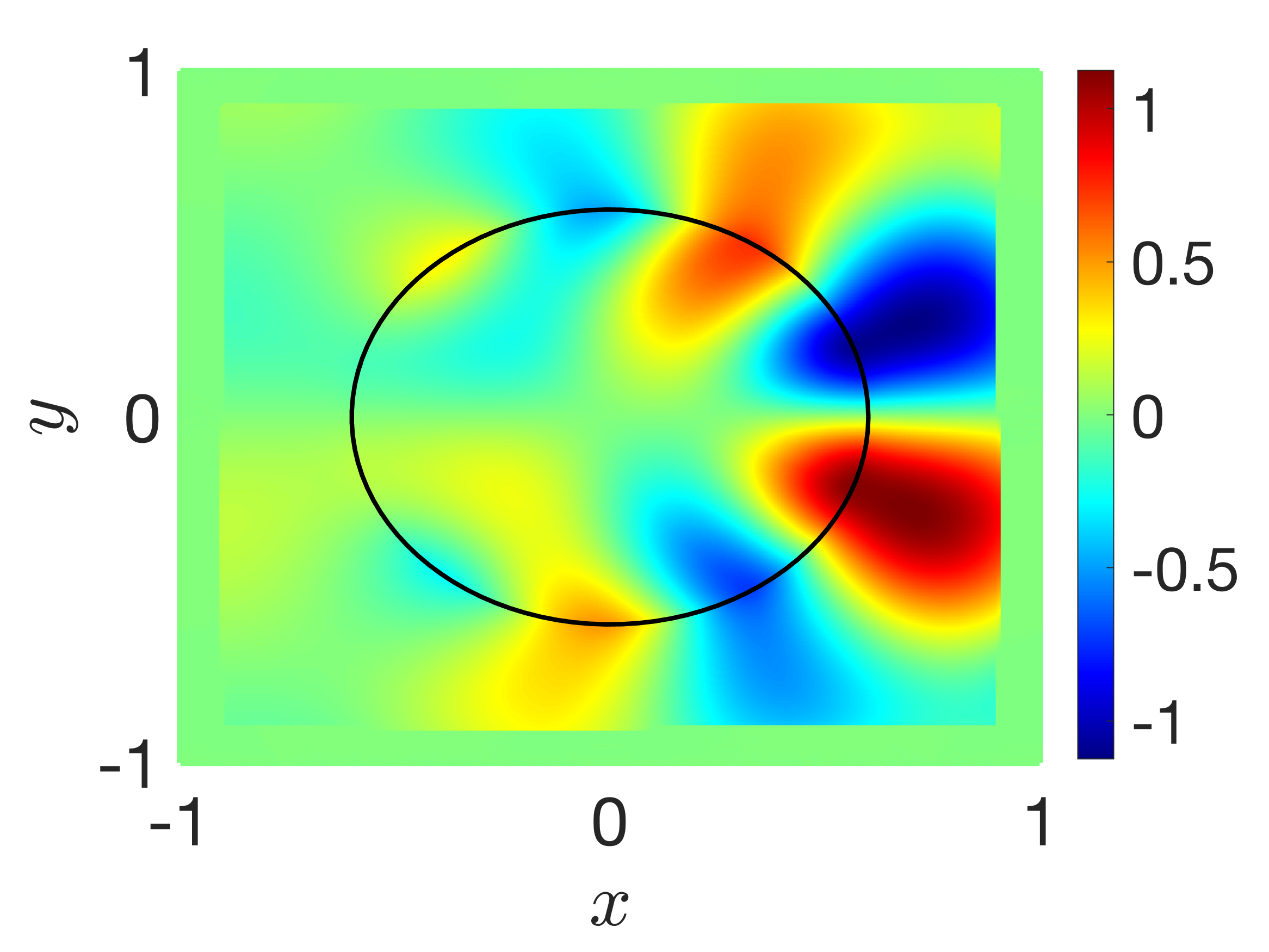}}{\footnotesize$H_x$}
	\stackunder[5pt]{	\includegraphics[width=1.58in]{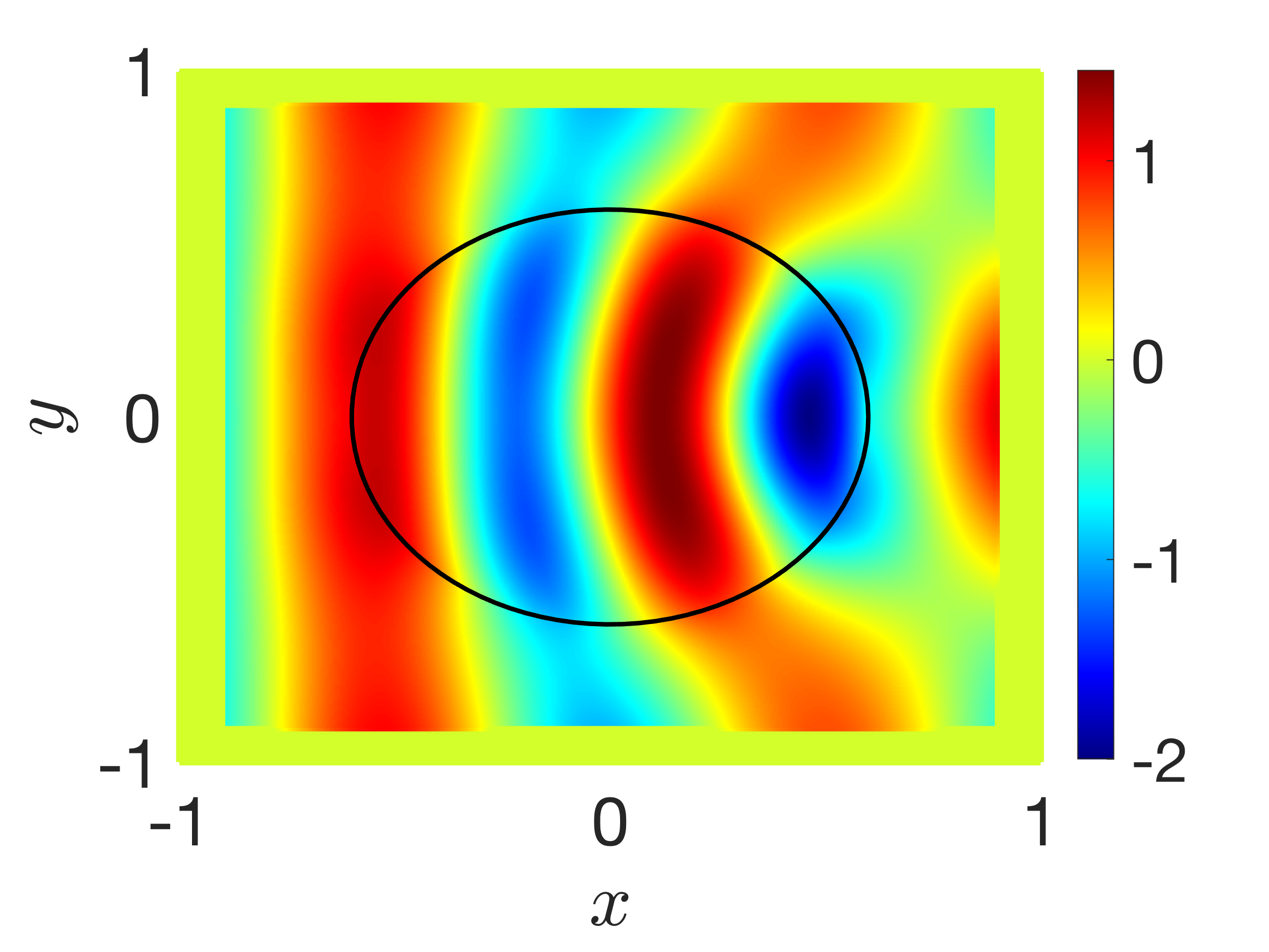}}{\footnotesize$H_y$}
	\stackunder[5pt]{	\includegraphics[width=1.58in]{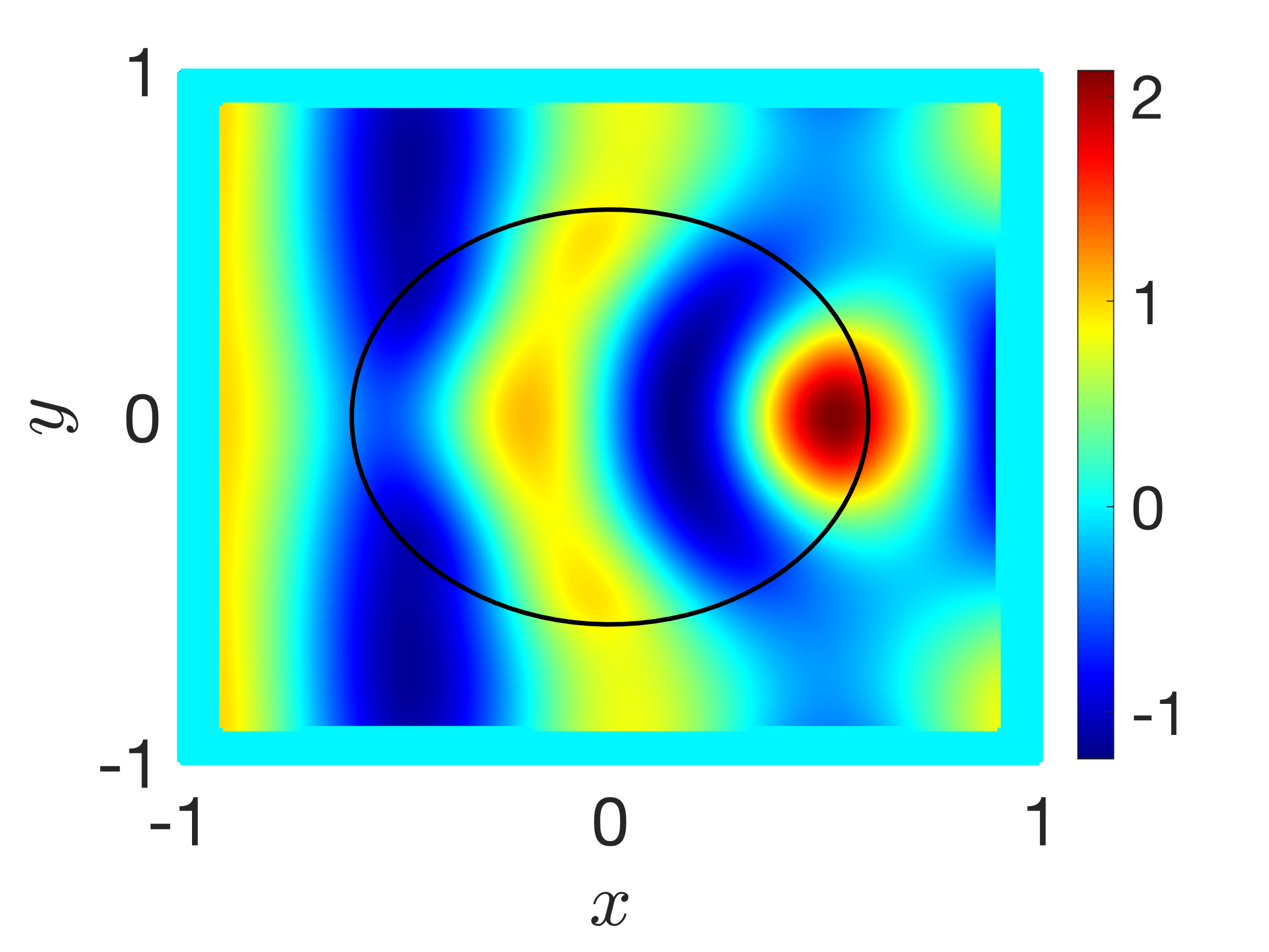}}{\footnotesize$E_z$}
	} 
	{\phantom{A}}
	\label{fig:nonMagneticSols}
	}
	\subfigure[a magnetic dielectric material ($\mu^+=1$, $\mu^-=2$, $\epsilon^+=1$ and $\epsilon^-=2.25$)]{
	\stackunder[5pt]{
	\stackunder[5pt]{\includegraphics[width=1.58in]{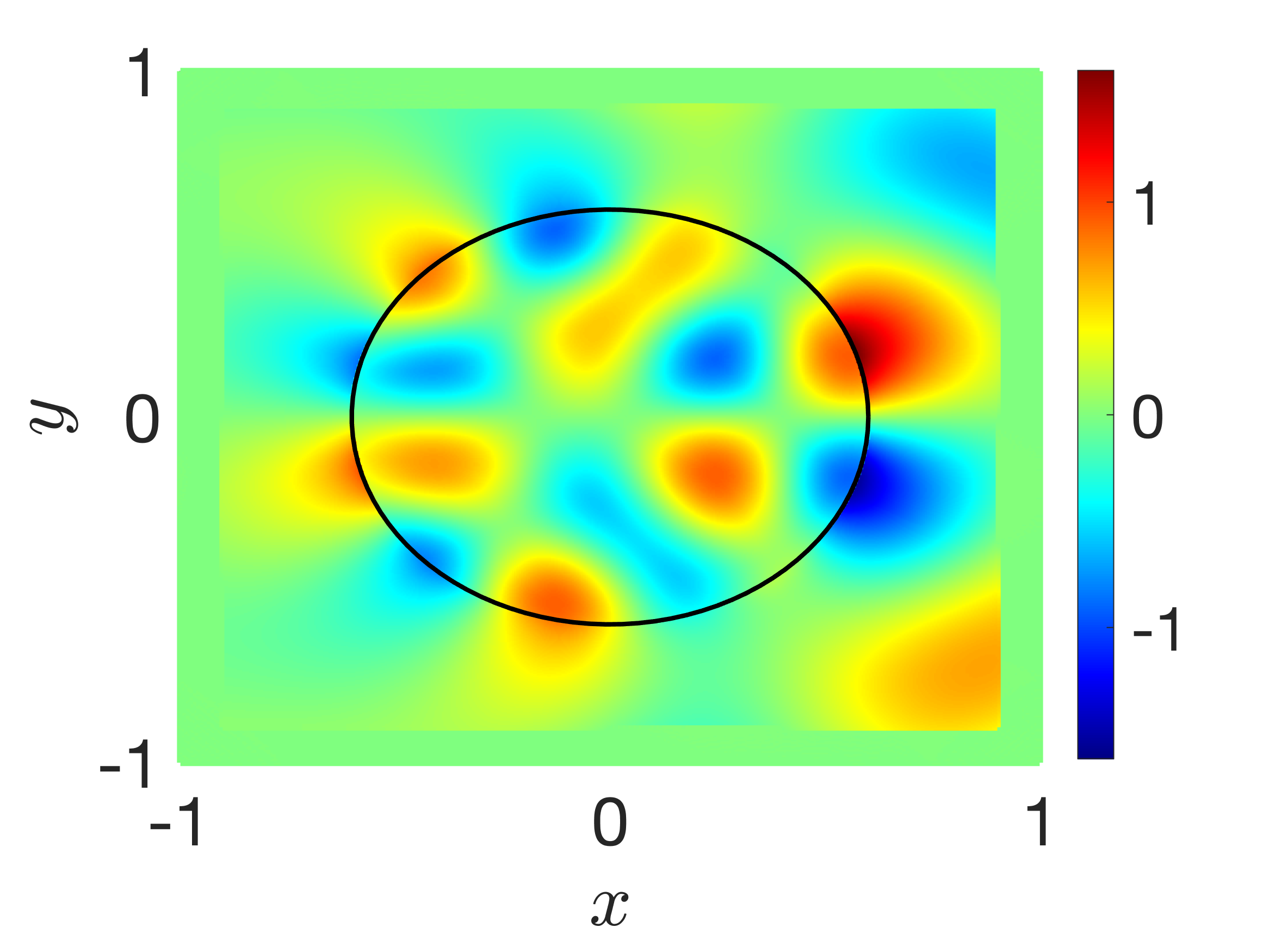}}{\footnotesize$H_x$}
	\stackunder[5pt]{	\includegraphics[width=1.58in]{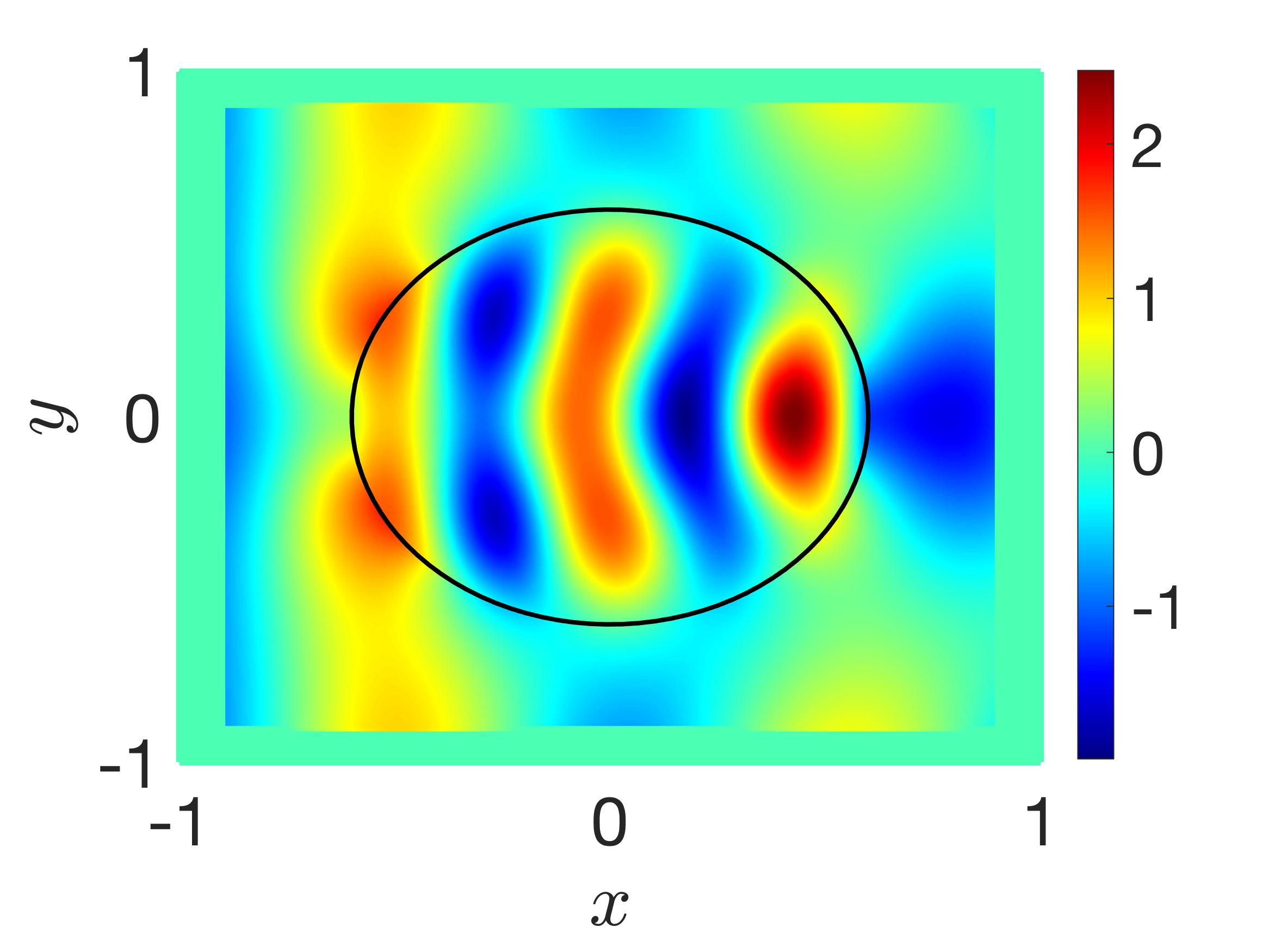}}{\footnotesize$H_y$}
	\stackunder[5pt]{	\includegraphics[width=1.58in]{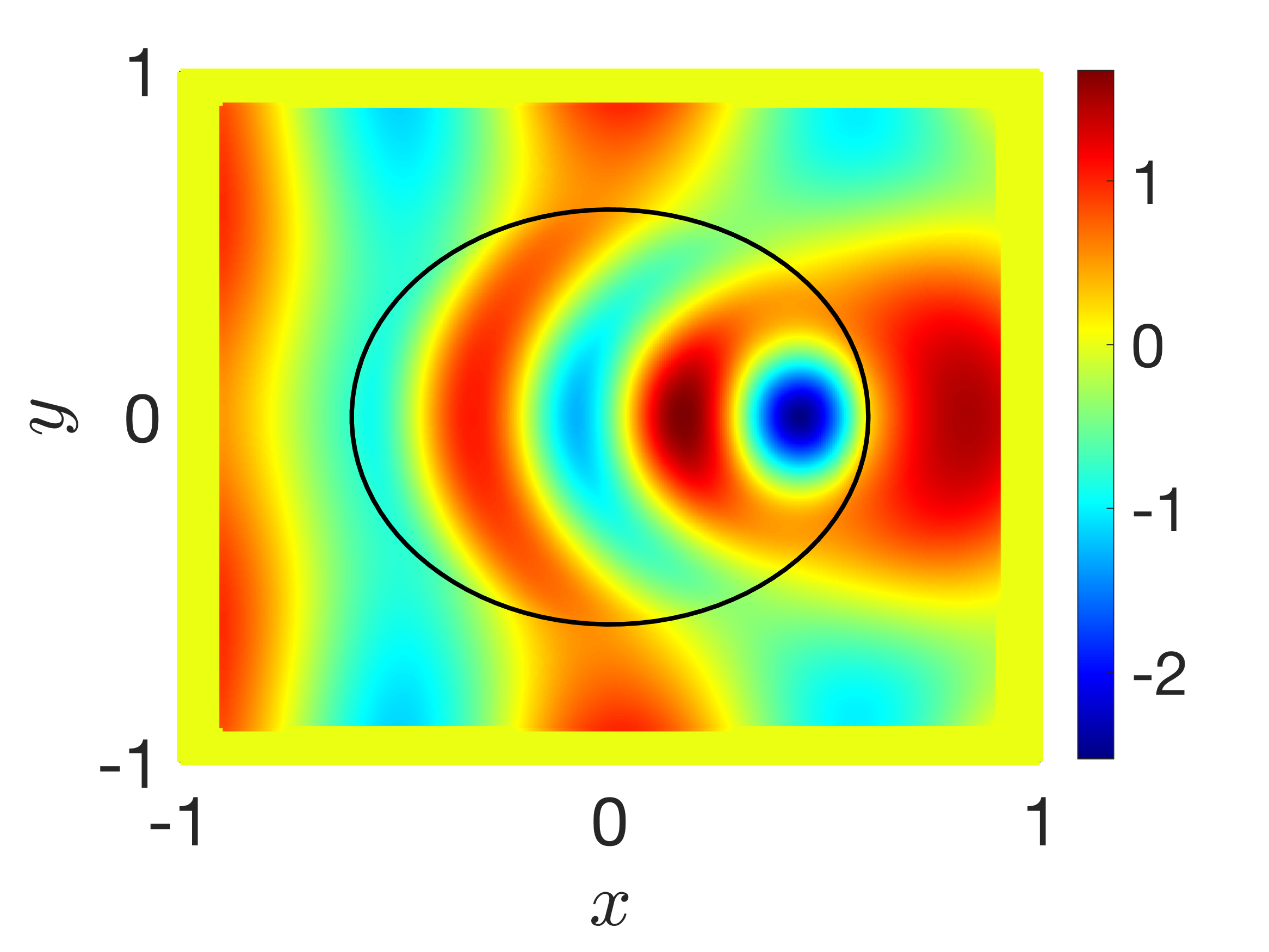}}{\footnotesize$E_z$}
	} 
	{\phantom{A}}
	\label{fig:magneticSols}
	}
       \caption{The components $H_x$, $H_y$ and $E_z$ with $h = \tfrac{1}{244}$ for scattering of a dielectric cylinder problems using the CFM-$4^{th}$ scheme. The computed electric field and 
       magnetic field are shown respectively at $t=1$ and $t-\tfrac{\Delta t}{2}$. The interface is represented by the black line.}
       \label{fig:plotInterfaceFieldsAnalyticSol}
\end{figure}

Let us first consider $\mu^+=\mu^-=1$, 
	$\epsilon^+=1$ and $\epsilon^-=2.25$. 
This corresponds to a non-magnetic dielectric material,
	and therefore $H_x$, 
	$H_y$ and $E_z$ are continuous across the interface. 
 \begin{figure} 
 \centering
  \subfigure[non-magnetic case ($\mu^-=1$)]{ \label{fig:convPlotNonMagnet}
		\setlength\figureheight{0.33\linewidth} 
		\setlength\figurewidth{0.33\linewidth} 
		\tikzset{external/export next=false}
\begin{tikzpicture}

\begin{axis}[%
width=0.951\figurewidth,
height=\figureheight,
at={(0\figurewidth,0\figureheight)},
scale only axis,
xmode=log,
xmin=0.00297619047619048,
xmax=0.1,
xminorticks=true,
xlabel style={font=\color{white!15!black}},
xlabel={\scriptsize$h$},
ymode=log,
ymin=1e-08,
ymax=1,
yminorticks=true,
ytick = {1e-8,1e-7,1e-6,1e-5, 1e-4 ,1e-3, 1e-2, 1e-1,1},
ylabel={\scriptsize$\|\mathbold{U}-\mathbold{U}_h\|_{2}$},
axis background/.style={fill=white},
legend style={at={(0.62,0.5)},anchor=north west,legend cell align=left,align=left,draw=white!15!black,draw=none,fill=none},
legend style={font=\scriptsize},
ylabel style={yshift=-5pt},xlabel style={yshift=2.5pt},tick label style={font=\tiny} 
]
\addplot [color=black,line width=1pt,solid,mark=o,mark options={solid}]
  table[row sep=crcr]{%
0.05	0.0597405133754351\\
0.0357142857142857	0.0315808675647812\\
0.025	0.0156232291006089\\
0.0192307692307692	0.00937734672857387\\
0.0138888888888889	0.00493031410509885\\
0.0104166666666667	0.0027861608747584\\
0.00757575757575758	0.00148724984074204\\
0.00555555555555556	0.000803005065359366\\
0.00409836065573771	0.000437209093036551\\
0.00297619047619048	0.000231604451276791\\
};
\addlegendentry{Yee}

\addplot [color=blue,line width=1pt,solid,mark=square,mark options={solid}]
  table[row sep=crcr]{%
0.05	0.0151965035860484\\
0.0357142857142857	0.0043325826882963\\
0.025	0.00115534011483271\\
0.0192307692307692	0.000412589460535316\\
0.0138888888888889	0.000116553241995696\\
0.0104166666666667	3.86518908289245e-05\\
0.00757575757575758	1.12088879438502e-05\\
0.00555555555555556	3.29937662687191e-06\\
0.00409836065573771	1.01449393963461e-06\\
0.002976190476190	2.794443112649193e-07\\
};
\addlegendentry{4$^{th}$}

\addplot [color=red,line width=1pt,densely dashed]
  table[row sep=crcr]{%
0.05	0.25\\
0.0357142857142857	0.127551020408163\\
0.025	0.0625\\
0.0192307692307692	0.0369822485207101\\
0.0138888888888889	0.0192901234567901\\
0.0104166666666667	0.0108506944444444\\
0.00757575757575758	0.00573921028466483\\
0.00555555555555556	0.00308641975308642\\
0.00409836065573771	0.00167965600644988\\
0.00297619047619048	0.000885770975056689\\
};
\addlegendentry{$h^2$}

\addplot [color=red,line width=1pt,solid]
  table[row sep=crcr]{%
0.05	0.005\\
0.0357142857142857	0.00130154102457309\\
0.025	0.0003125\\
0.0192307692307692	0.000109414936451805\\
0.0138888888888889	2.97687090382564e-05\\
0.0104166666666667	9.4190055941358e-06\\
0.00757575757575758	2.63508277532821e-06\\
0.00555555555555556	7.62078951379363e-07\\
0.00409836065573771	2.25699544000253e-07\\
0.00297619047619048	6.27672176202302e-08\\
};
\addlegendentry{$h^4$}

\end{axis}

\end{tikzpicture}%
		}		
  \subfigure[magnetic case ($\mu^-=2$)]{ \label{fig:convPlotMagnet}
		\setlength\figureheight{0.33\linewidth} 
		\setlength\figurewidth{0.33\linewidth} 
		\tikzset{external/export next=false}
\begin{tikzpicture}

\begin{axis}[%
width=0.951\figurewidth,
height=\figureheight,
at={(0\figurewidth,0\figureheight)},
scale only axis,
xmode=log,
xmin=0.00297619047619048,
xmax=0.1,
xminorticks=true,
xlabel style={font=\color{white!15!black}},
xlabel={\scriptsize$h$},
ymode=log,
ymin=1e-07,
ymax=1,
yminorticks=true,
ytick = {1e-8,1e-7,1e-6,1e-5, 1e-4 ,1e-3, 1e-2, 1e-1,1},
ylabel={\scriptsize$\|\mathbold{U}-\mathbold{U}_h\|_{2}$},
axis background/.style={fill=white},
legend style={at={(0.62,0.5)},anchor=north west,legend cell align=left,align=left,draw=white!15!black,draw=none,fill=none},
legend style={font=\scriptsize},
ylabel style={yshift=-5pt},xlabel style={yshift=2.5pt},tick label style={font=\tiny} 
]
\addplot [color=black,line width=1pt,solid,mark=o,mark options={solid}]
  table[row sep=crcr]{%
0.05	0.122639202549746\\
0.0357142857142857	0.0632486022174324\\
0.025	0.0311470892042609\\
0.0192307692307692	0.0186419717220955\\
0.0138888888888889	0.00976201541030907\\
0.0104166666666667	0.00552557135163722\\
0.00757575757575758	0.00294280237319364\\
0.00555555555555556	0.00158461435475192\\
0.00409836065573771	0.00086403207712408\\
0.00297619047619048	0.000456870226163026\\
};
\addlegendentry{Yee}

\addplot [color=blue,line width=1pt,solid,mark=square,mark options={solid}]
  table[row sep=crcr]{%
0.05	0.0174790260866363\\
0.0357142857142857	0.00648199611770236\\
0.025	0.00178436385182201\\
0.0192307692307692	0.000628850133049453\\
0.0138888888888889	0.000209104676139261\\
0.0104166666666667	6.69393493959338e-05\\
0.00757575757575758	1.89597629583939e-05\\
0.00555555555555556	6.04094930829162e-06\\
0.00409836065573771	1.93319610702559e-06\\
0.002976190476190	        4.809616668554840e-07\\
};
\addlegendentry{4$^{th}$}

\addplot [color=red,line width=1pt,densely dashed]
  table[row sep=crcr]{%
0.05	0.4375\\
0.0357142857142857	0.223214285714286\\
0.025	0.109375\\
0.0192307692307692	0.0647189349112426\\
0.0138888888888889	0.0337577160493827\\
0.0104166666666667	0.0189887152777778\\
0.00757575757575758	0.0100436179981635\\
0.00555555555555556	0.00540123456790123\\
0.00409836065573771	0.00293939801128729\\
0.00297619047619048	0.00155009920634921\\
};
\addlegendentry{$h^2$}

\addplot [color=red,line width=1pt,solid]
  table[row sep=crcr]{%
0.05	0.009375\\
0.0357142857142857	0.00244038942107455\\
0.025	0.0005859375\\
0.0192307692307692	0.000205153005847134\\
0.0138888888888889	5.58163294467307e-05\\
0.0104166666666667	1.76606354890046e-05\\
0.00757575757575758	4.94078020374038e-06\\
0.00555555555555556	1.42889803383631e-06\\
0.00409836065573771	4.23186645000474e-07\\
0.00297619047619048	1.17688533037932e-07\\
};
\addlegendentry{$h^4$}

\end{axis}
\end{tikzpicture}%
  }
  \caption{Convergence plots for scattering of a dielectric cylinder problems with $\mu^+=1$, $\epsilon^+=1$ and $\epsilon^-=2.25$ using the proposed CFM-FDTD schemes. It is recalled that $\mathbold{U} = [H_x,H_y,E_z]^T$.}
   \label{fig:convPlotScatteringDielectricPblms}
\end{figure}
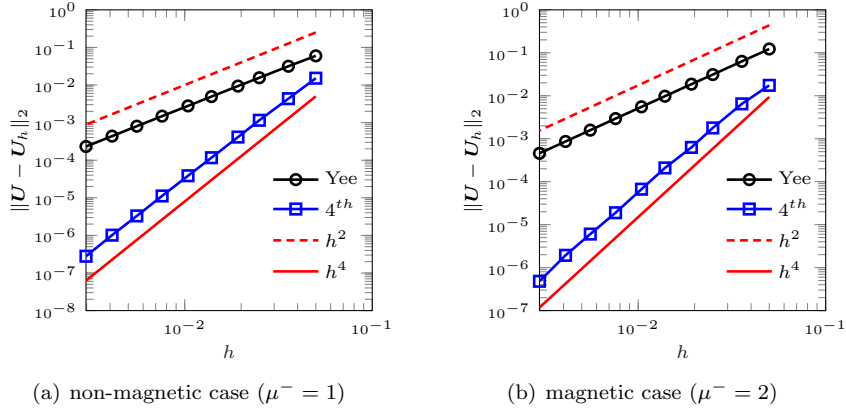
Figure~\ref{fig:convPlotNonMagnet} illustrates the convergence plot of 
	$\mathbold{U} = [H_x,H_y,E_z]^T$ for both CFM-FDTD schemes.
We observe a second-order convergence in $L^2$-norm for the CFM-Yee scheme as expected by the theory.
For the CFM-4$^{th}$ scheme, 
	a fourth-order convergence is obtained, 
	which is better than expected.
Numerical solutions computed with the CFM-4$^{th}$ scheme at $t=1$ are illustrated 
	in Figure~\ref{fig:nonMagneticSols}.
	
Let us now consider a magnetic dielectric material. 
We choose $\mu^+ = 1$,
	$\mu^-=2$,
	$\epsilon^+ = 1$ and $\epsilon^-=2.25$. 
In this case, 
	the components of the magnetic field are discontinuous 
	while the $z$-component of the electric field is still continuous across the interface.
Figure~\ref{fig:convPlotMagnet} illustrates the convergence plot of electromagnetic fields
	for both schemes.	
A second and fourth order convergence in $L^2$-norm are observed for respectively the CFM-Yee and 
	the CFM-$4^{th}$ scheme.
These results are in agreement with the theory.
Figure~\ref{fig:magneticSols} illustrates the approximation of $H_x$, 
	$H_y$ and $E_z$ at $t=1$.

%***************************************
%***************************************
	
\subsubsection{Verification of the Accuracy of Correction Functions}
In this subsection, 
	we assess the accuracy of the estimated correction functions coming from 
	minimization problem \eqref{eq:minPblm} using high-order explicit jump 
	conditions \cite{Zhao2004,Zhang2016}. 
Matched Interface and Boundary (MIB) based strategies 
	use these conditions to construct high-order FDTD schemes.
As mentioned in Remark~\ref{rem:jumpCdns}, 
	the correction function's system of PDEs implicitly considers 
	jump conditions coming from Maxwell's equations \eqref{eq:pblmDefinition}. 
To provide numerical evidences of this claim, 
	we compute the error on these jump conditions on all local patches using 
\begin{equation*}
	\Big(\int\limits_{\Gamma\cap\Omega_{\Gamma}^{h}}\!\!\llbracket u(\mathbold{x},t_f) \rrbracket^2 \, \mathrm{d}S\Big)^{1/2},
\end{equation*}
	where $\llbracket u(\mathbold{x},t_f) \rrbracket$ is a given jump condition evaluated with 
	approximated solutions of problem~\eqref{eq:minPblm} at $t_f$. 
Afterward, 
	the maximum error value on all local patches for a given order of jump 
	conditions is taken and is denoted by $E_i$ for the $i^{th}$-order jump condition.
	
Although we do not have a theoretical result to characterize the convergence of 
	high-order explicit jump conditions, 
	one should expect a $(k+1-q)$ convergence for a $q^{th}$ order jump condition when  
	$k$ degree polynomial approximations of correction functions are used.
As an example, 
	a third degree polynomial approximation should lead at least to a fourth, 
	third, 
	second 
	and first order convergence for respectively the zeroth, 
	first,
	second 
	and third order jump conditions.
It is recalled that second and third degree polynomial approximations of correction functions 
	are used for respectively the CFM-Yee scheme and the CFM-$4^{th}$ scheme.
	
\begin{figure} 
 \centering
  \subfigure[zeroth-order jump conditions]{ \label{fig:zerothOrder}
		\setlength\figureheight{0.33\linewidth} 
		\setlength\figurewidth{0.33\linewidth} 
		\tikzset{external/export next=false}
\begin{tikzpicture}

\begin{axis}[%
width=0.951\figurewidth,
height=\figureheight,
at={(0\figurewidth,0\figureheight)},
scale only axis,
xmode=log,
xmin=0.00297619047619048,
xmax=0.1,
xminorticks=true,
xlabel style={font=\color{white!15!black}},
xlabel={\scriptsize$h$},
ymode=log,
ymin=1e-09,
ymax=0.01,
yminorticks=true,
ytick = {1e-9,1e-8,1e-7,1e-6,1e-5, 1e-4 ,1e-3, 1e-2, 1e-1,1},
ylabel={\scriptsize$E_0$},
axis background/.style={fill=white},
legend style={at={(0.615,0.5)},anchor=north west,legend cell align=left,align=left,draw=white!15!black,draw=none,fill=none},
legend style={font=\scriptsize},
ylabel style={yshift=-5pt},xlabel style={yshift=2.5pt},tick label style={font=\tiny} 
]
\addplot [color=black,line width=1pt,solid,mark=o,mark options={solid}]
  table[row sep=crcr]{%
0.05	0.00337362976965725\\
0.0357142857142857	0.000801307015503004\\
0.025	0.000229065122091851\\
0.0192307692307692	8.96471357098544e-05\\
0.0138888888888889	2.89198687082825e-05\\
0.0104166666666667	1.06242355089369e-05\\
0.00757575757575758	3.47734375026281e-06\\
0.00555555555555556	1.15901841359661e-06\\
0.00409836065573771	3.98857978501739e-07\\
0.00297619047619048	1.29630866252407e-07\\
};
\addlegendentry{Yee}

\addplot [color=blue,line width=1pt,solid,mark=square,mark options={solid}]
  table[row sep=crcr]{%
0.05	0.00181910312788471\\
0.0357142857142857	0.000476112706850772\\
0.025	0.000105925847406062\\
0.0192307692307692	3.41340014281073e-05\\
0.0138888888888889	8.19871611164324e-06\\
0.0104166666666667	2.29559745348512e-06\\
0.00757575757575758	5.58376909517785e-07\\
0.00555555555555556	1.39891461432322e-07\\
0.00409836065573771	3.59250319148011e-08\\
0.00297619047619048	8.56496952410249e-09\\
};
\addlegendentry{4$^{th}$}

\addplot [color=red,line width=1pt,densely dashed]
  table[row sep=crcr]{%
0.05	0.00838525491562421\\
0.0357142857142857	0.00258266540412646\\
0.025	0.000741158826601964\\
0.0192307692307692	0.000295875872758962\\
0.0138888888888889	9.47234521265415e-05\\
0.0104166666666667	3.46076363850795e-05\\
0.00757575757575758	1.13530485612646e-05\\
0.00555555555555556	3.83413576388853e-06\\
0.00409836065573771	1.32207734278519e-06\\
0.00297619047619048	4.31453288859707e-07\\
};
\addlegendentry{$h^{3.5}$}

\addplot [color=red,line width=1pt,solid]
  table[row sep=crcr]{%
0.05	0.000698771242968685\\
0.0357142857142857	0.000153730083578956\\
0.025	3.08816177750818e-05\\
0.0192307692307692	9.48320104996674e-06\\
0.0138888888888889	2.1926725029292e-06\\
0.0104166666666667	6.00827020574297e-07\\
0.00757575757575758	1.4334657274324e-07\\
0.00555555555555556	3.55012570730419e-08\\
0.00409836065573771	9.03058294252178e-09\\
0.00297619047619048	2.1401452820422e-09\\
};
\addlegendentry{$h^{4.5}$}

\end{axis}

\end{tikzpicture}%
  }
  \subfigure[first-order jump conditions]{ \label{fig:firstOrder}
		\setlength\figureheight{0.33\linewidth} 
		\setlength\figurewidth{0.33\linewidth} 
		\tikzset{external/export next=false}
\begin{tikzpicture}

\begin{axis}[%
width=0.951\figurewidth,
height=\figureheight,
at={(0\figurewidth,0\figureheight)},
scale only axis,
xmode=log,
xmin=0.00297619047619048,
xmax=0.1,
xminorticks=true,
xlabel style={font=\color{white!15!black}},
xlabel={\scriptsize$h$},
ymode=log,
ymin=1e-06,
ymax=10,
yminorticks=true,
ytick = {1e-8,1e-7,1e-6,1e-5, 1e-4 ,1e-3, 1e-2, 1e-1,1},
ylabel={\scriptsize$E_1$},
axis background/.style={fill=white},
legend style={at={(0.615,0.5)},anchor=north west,legend cell align=left,align=left,draw=white!15!black,draw=none,fill=none},
legend style={font=\scriptsize},
ylabel style={yshift=-5pt},xlabel style={yshift=2.5pt},tick label style={font=\tiny} 
]
\addplot [color=black,line width=1pt,solid,mark=o,mark options={solid}]
  table[row sep=crcr]{%
0.05	0.33220182483918\\
0.0357142857142857	0.161266908081682\\
0.025	0.0685464727366656\\
0.0192307692307692	0.0296992742359436\\
0.0138888888888889	0.0128252644286622\\
0.0104166666666667	0.00562958644676737\\
0.00757575757575758	0.00226488259008506\\
0.00555555555555556	0.00092537291053364\\
0.00409836065573771	0.000384782735312521\\
0.00297619047619048	0.000139227285425775\\
};
\addlegendentry{Yee}

\addplot [color=blue,line width=1pt,solid,mark=square,mark options={solid}]
  table[row sep=crcr]{%
0.05	0.115391081562178\\
0.0357142857142857	0.0408040613874252\\
0.025	0.0104122184458359\\
0.0192307692307692	0.00444272344565276\\
0.0138888888888889	0.00136600841328604\\
0.0104166666666667	0.000465118167378365\\
0.00757575757575758	0.00015573552704794\\
0.00555555555555556	5.18570897911786e-05\\
0.00409836065573771	1.76787641571236e-05\\
0.00297619047619048	5.72696514054515e-06\\
};
\addlegendentry{4$^{th}$}

\addplot [color=red,line width=1pt,densely dashed]
  table[row sep=crcr]{%
0.05	1.625\\
0.0357142857142857	0.592201166180758\\
0.025	0.203125\\
0.0192307692307692	0.0924556213017752\\
0.0138888888888889	0.0348293895747599\\
0.0104166666666667	0.0146936487268518\\
0.00757575757575758	0.005652252553079\\
0.00555555555555556	0.00222908093278464\\
0.00409836065573771	0.000894898691961001\\
0.00297619047619048	0.000342709008206457\\
};
\addlegendentry{$h^3$}

\addplot [color=red,line width=1pt,solid]
  table[row sep=crcr]{%
0.05	0.0363361046343716\\
0.0357142857142857	0.011191550084548\\
0.025	0.00321168824860851\\
0.0192307692307692	0.0012821287819555\\
0.0138888888888889	0.000410468292548347\\
0.0104166666666667	0.000149966424335345\\
0.00757575757575758	4.919654376548e-05\\
0.00555555555555556	1.66145883101836e-05\\
0.00409836065573771	5.72900181873582e-06\\
0.00297619047619048	1.86963091839206e-06\\
};
\addlegendentry{$h^{3.5}$}

\end{axis}
\end{tikzpicture}%
		}		
  \subfigure[second-order jump conditions]{ \label{fig:secondOrder}
		\setlength\figureheight{0.33\linewidth} 
		\setlength\figurewidth{0.33\linewidth} 
		\tikzset{external/export next=false}
\begin{tikzpicture}

\begin{axis}[%
width=0.951\figurewidth,
height=\figureheight,
at={(0\figurewidth,0\figureheight)},
scale only axis,
xmode=log,
xmin=0.00297619047619048,
xmax=0.1,
xminorticks=true,
xlabel style={font=\color{white!15!black}},
xlabel={\scriptsize$h$},
ymode=log,
ymin=0.0001,
ymax=100,
yminorticks=true,
ytick = {1e-8,1e-7,1e-6,1e-5, 1e-4 ,1e-3, 1e-2, 1e-1,1,10,100,1000},
ylabel={\scriptsize$E_2$},
axis background/.style={fill=white},
legend style={at={(0.615,0.5)},anchor=north west,legend cell align=left,align=left,draw=white!15!black,draw=none,fill=none},
legend style={font=\scriptsize},
ylabel style={yshift=-5pt},xlabel style={yshift=2.5pt},tick label style={font=\tiny} 
]
\addplot [color=black,line width=1pt,solid,mark=o,mark options={solid}]
  table[row sep=crcr]{%
0.05	17.8215225904192\\
0.0357142857142857	11.7246463744881\\
0.025	7.13904457229311\\
0.0192307692307692	4.1665231551818\\
0.0138888888888889	2.41481808689406\\
0.0104166666666667	1.47167267146539\\
0.00757575757575758	0.785270589858474\\
0.00555555555555556	0.441606638183393\\
0.00409836065573771	0.254015772022084\\
0.00297619047619048	0.124595802290699\\
};
\addlegendentry{Yee}

\addplot [color=blue,line width=1pt,solid,mark=square,mark options={solid}]
  table[row sep=crcr]{%
0.05	5.60025244009404\\
0.0357142857142857	3.02090084788958\\
0.025	1.16509599933616\\
0.0192307692307692	0.520245315047039\\
0.0138888888888889	0.194778064243259\\
0.0104166666666667	0.0907034265749897\\
0.00757575757575758	0.0328145950989696\\
0.00555555555555556	0.012586273904851\\
0.00409836065573771	0.0053391942008075\\
0.00297619047619048	0.00241211245390837\\
};
\addlegendentry{4$^{th}$}

\addplot [color=red,line width=1pt,densely dashed]
  table[row sep=crcr]{%
0.05	75\\
0.0357142857142857	38.265306122449\\
0.025	18.75\\
0.0192307692307692	11.094674556213\\
0.0138888888888889	5.78703703703704\\
0.0104166666666667	3.25520833333333\\
0.00757575757575758	1.72176308539945\\
0.00555555555555556	0.925925925925926\\
0.00409836065573771	0.503896801934964\\
0.00297619047619048	0.265731292517007\\
};
\addlegendentry{$h^2$}

\addplot [color=red,line width=1pt,solid]
  table[row sep=crcr]{%
0.05	4.375\\
0.0357142857142857	1.59438775510204\\
0.025	0.546875\\
0.0192307692307692	0.248918980427856\\
0.0138888888888889	0.0937714334705075\\
0.0104166666666667	0.0395598234953704\\
0.00757575757575758	0.0152176030275204\\
0.00555555555555556	0.00600137174211248\\
0.00409836065573771	0.0024093426322027\\
0.00297619047619048	0.000922678099017385\\
};
\addlegendentry{$h^3$}

\end{axis}

\end{tikzpicture}%
  }
  \subfigure[third-order jump conditions]{ \label{fig:thirdOrder}
		\setlength\figureheight{0.33\linewidth} 
		\setlength\figurewidth{0.33\linewidth} 
		\tikzset{external/export next=false}
\begin{tikzpicture}

\begin{axis}[%
width=0.951\figurewidth,
height=\figureheight,
at={(0\figurewidth,0\figureheight)},
scale only axis,
xmode=log,
xmin=0.00297619047619048,
xmax=0.1,
xminorticks=true,
xlabel style={font=\color{white!15!black}},
xlabel={\scriptsize$h$},
ymode=log,
ymin=0.442885487528345,
ymax=125,
yminorticks=true,
ytick = {1e-8,1e-7,1e-6,1e-5, 1e-4 ,1e-3, 1e-2, 1e-1,1,10,100,1000},
ylabel={\scriptsize$E_3$},
axis background/.style={fill=white},
legend style={at={(0.615,0.5)},anchor=north west,legend cell align=left,align=left,draw=white!15!black,draw=none,fill=none},
legend style={font=\scriptsize},
ylabel style={yshift=-5pt},xlabel style={yshift=2.5pt},tick label style={font=\tiny} 
]
\addplot [color=blue,line width=1pt,solid,mark=square,mark options={solid}]
  table[row sep=crcr]{%
0.05	114.667051177018\\
0.0357142857142857	92.4751838764351\\
0.025	54.5451991828081\\
0.0192307692307692	31.9614523781774\\
0.0138888888888889	16.145694431744\\
0.0104166666666667	10.1884632992507\\
0.00757575757575758	4.95590419914097\\
0.00555555555555556	2.58356468532146\\
0.00409836065573771	1.35555786866128\\
0.00297619047619048	0.966865333579854\\
};
\addlegendentry{4$^{th}$}

\addplot [color=red,line width=1pt,solid]
  table[row sep=crcr]{%
0.05	125\\
0.0357142857142857	63.7755102040816\\
0.025	31.25\\
0.0192307692307692	18.491124260355\\
0.0138888888888889	9.64506172839506\\
0.0104166666666667	5.42534722222222\\
0.00757575757575758	2.86960514233242\\
0.00555555555555556	1.54320987654321\\
0.00409836065573771	0.83982800322494\\
0.00297619047619048	0.442885487528345\\
};
\addlegendentry{$h^2$}

\end{axis}
\end{tikzpicture}%
  }
  \caption{Convergence plots of jump conditions for a scattering of a non-magnetic dielectric cylinder problem ($\mu^+=\mu^-=1$) using the proposed CFM-FDTD schemes.}
   \label{fig:convPlotNonMagneticJumpCdns}
\end{figure}
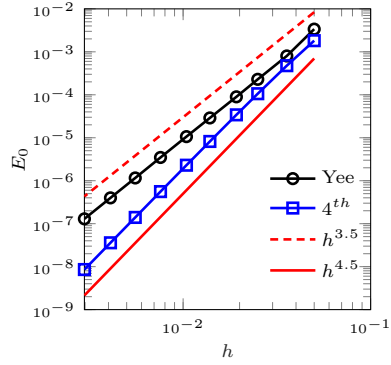
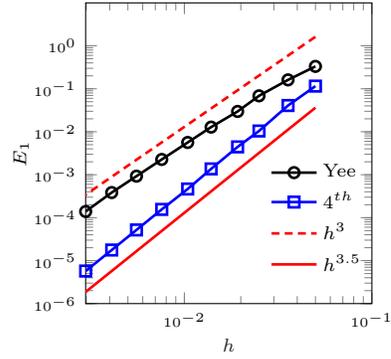
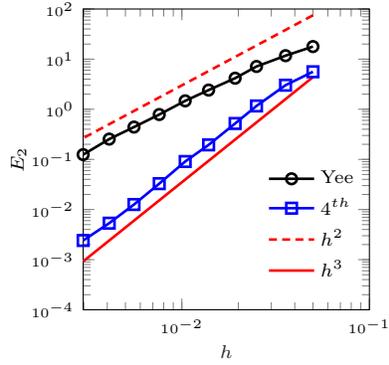
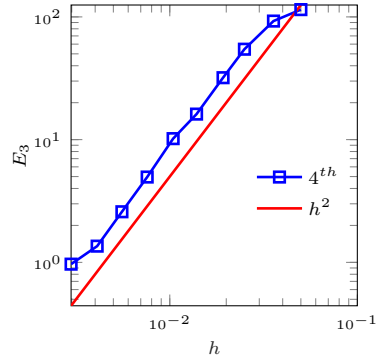

For a non-magnetic dielectric material, 
	high-order jump conditions can be derived by using the continuity of 
	time derivatives of electromagnetic fields on the interface \cite{Zhao2004} and  
	are given by: 
\begin{equation*}
\begin{aligned}
\text{zeroth-order} &\,\, \left\{ 
  \begin{array}{l}
    \llbracket H_x \rrbracket =0, \\ %[1.5mm]
    \llbracket H_y \rrbracket =0, \\
    \llbracket E_z \rrbracket =0, \\
  \end{array} \right. \qquad
\text{first-order} \,\, \left\{ 
  \begin{array}{l l}
    \llbracket \partial_y E_z \rrbracket =0, \\ %[1.5mm]
    \llbracket \partial_x E_z \rrbracket =0, \\
    \big\llbracket \tfrac{1}{\epsilon}\,(\partial_xH_y-\partial_yH_x) \big\rrbracket =0, \\
  \end{array} \right. \\
\text{second-order} &\,\, \left\{ 
  \begin{array}{l l}
    \big\llbracket \tfrac{1}{\epsilon}\,(\partial_x^2E_z-\partial_y^2E_z) \big\rrbracket =0, \\ %[1.5mm]
   \big\llbracket \tfrac{1}{\epsilon}\,(\partial_y^2H_x-\partial_{xy}^2H_y) \big\rrbracket =0,\\
   \big\llbracket \tfrac{1}{\epsilon}\,(\partial_x^2H_y-\partial_{xy}^2H_x) \big\rrbracket =0,\\
  \end{array} \right. \\
\text{third-order} &\,\, \left\{ 
  \begin{array}{l l}
    \big\llbracket \tfrac{1}{\epsilon}\,(\partial_{xxy}^3E_z+\partial_y^3E_z) \big\rrbracket =0, \\ %[1.5mm]
   \big\llbracket \tfrac{1}{\epsilon}\,(\partial_x^3E_z+\partial_{xyy}^3E_z) \big\rrbracket =0,\\
   \big\llbracket \tfrac{1}{\epsilon^2}\,(\partial_x^3H_y+\partial_{xyy}^3H_y-\partial_y^3H_x-\partial_{xxy}^3H_x) \big\rrbracket =0.
  \end{array} \right.
  \end{aligned}
\end{equation*}
Figure~\ref{fig:convPlotNonMagneticJumpCdns} illustrates convergence plots of 
	those jump conditions at $t_f = 1$ for both schemes.
We observe that the convergence order for all jump conditions is better than expected. 

Let us now consider a magnetic dielectric material. 
Considering a point $\mathbold{p}=(x_p,y_p)$ on the interface $\Gamma$, 
	one can define a local coordinate system based on the normal $\mathbold{n}$ 
	and the tangent $\mathbold{\tau}$ to the interface at $\mathbold{p}$,
	and derive explicit jump conditions coming from Maxwell's 
	equations \eqref{eq:pblmDefinition} \cite{Zhang2016}.
In this local coordinate system, 
	zeroth and first order jump conditions are given by
\begin{equation*}
\begin{aligned}
\text{zeroth-order} &\,\, \left\{ 
  \begin{array}{l}
    \llbracket H_\tau \rrbracket =0, \\ %[1.5mm]
    \llbracket \mu\,H_n \rrbracket =0, \\
    \llbracket E_z \rrbracket =0, \\
  \end{array} \right. \\ %\qquad
\text{first-order} &\,\, \left\{ 
  \begin{array}{l l}
    \llbracket \partial_\tau E_z \rrbracket =0, \\ %[1.5mm]
    \llbracket \tfrac{1}{\mu}\,\partial_n E_z \rrbracket =0, \\
    \big\llbracket  \partial_n (\mu\,H_n)  + \partial_\tau (\mu\,H_\tau) \big\rrbracket =0, \\
    \big\llbracket \partial_n (\mu\,H_\tau)  -  \partial_\tau (\mu\,H_n)  -  \partial_t (\mu\,\epsilon\,E_z) \big\rrbracket=0.
%    \big\llbracket \partial_\tau (\mu\,H_\tau) \big\rrbracket + \big\llbracket \partial_n (\mu\,H_n) \big\rrbracket =0, \\
%    \big\llbracket \partial_n (\mu\,H_\tau) \big\rrbracket - \big\llbracket \partial_\tau (\mu\,H_n) \big\rrbracket - \big\llbracket \partial_t (\mu\,\epsilon\,E_z) \big\rrbracket=0.
  \end{array} \right. 
  \end{aligned}
\end{equation*}
 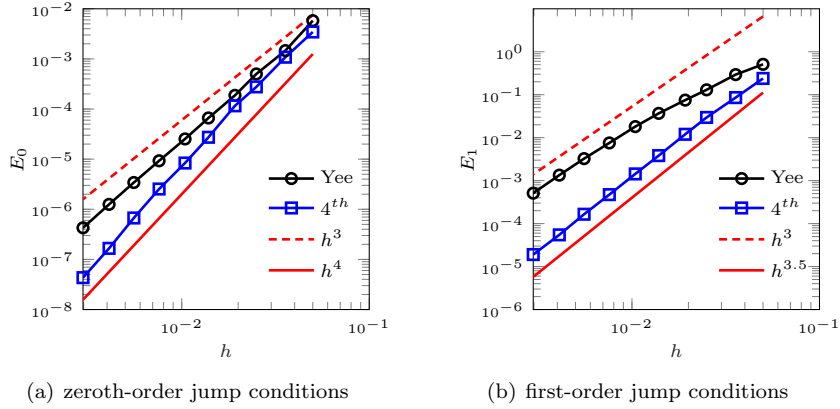
\begin{figure} 
 \centering
  \subfigure[zeroth-order jump conditions]{ \label{fig:zerothOrder}
		\setlength\figureheight{0.33\linewidth} 
		\setlength\figurewidth{0.33\linewidth} 
		\tikzset{external/export next=false}
\begin{tikzpicture}

\begin{axis}[%
width=0.951\figurewidth,
height=\figureheight,
at={(0\figurewidth,0\figureheight)},
scale only axis,
xmode=log,
xmin=0.00297619047619048,
xmax=0.1,
xminorticks=true,
xlabel style={font=\color{white!15!black}},
xlabel={\scriptsize$h$},
ymode=log,
ymin=1e-08,
ymax=0.01,
yminorticks=true,
ytick = {1e-8,1e-7,1e-6,1e-5, 1e-4 ,1e-3, 1e-2, 1e-1,1},
ylabel={\scriptsize$E_0$},
axis background/.style={fill=white},
legend style={at={(0.615,0.5)},anchor=north west,legend cell align=left,align=left,draw=white!15!black,draw=none,fill=none},
legend style={font=\scriptsize},
ylabel style={yshift=-5pt},xlabel style={yshift=2.5pt},tick label style={font=\tiny} 
]
\addplot [color=black,line width=1pt,solid,mark=o,mark options={solid}]
  table[row sep=crcr]{%
0.05	0.0057937290816341\\
0.0357142857142857	0.00146527709513542\\
0.025	0.000501066540117971\\
0.0192307692307692	0.000187623806964121\\
0.0138888888888889	6.63395582808436e-05\\
0.0104166666666667	2.53078197272906e-05\\
0.00757575757575758	9.27592075638467e-06\\
0.00555555555555556	3.42082964143953e-06\\
0.00409836065573771	1.25299979317757e-06\\
0.00297619047619048	4.26745417589557e-07\\
};
\addlegendentry{Yee}

\addplot [color=blue,line width=1pt,solid,mark=square,mark options={solid}]
  table[row sep=crcr]{%
0.05	0.00344770999844857\\
0.0357142857142857	0.0010823298144231\\
0.025	0.000276521910108713\\
0.0192307692307692	0.000116023609628376\\
0.0138888888888889	2.73300883916749e-05\\
0.0104166666666667	8.37450336208024e-06\\
0.00757575757575758	2.53777705117791e-06\\
0.00555555555555556	6.70889999046019e-07\\
0.00409836065573771	1.66229046049443e-07\\
0.00297619047619048	4.32617952367636e-08\\
};
\addlegendentry{4$^{th}$}

\addplot [color=red,line width=1pt,densely dashed]
  table[row sep=crcr]{%
0.05	0.0075\\
0.0357142857142857	0.0027332361516035\\
0.025	0.0009375\\
0.0192307692307692	0.000426718252162039\\
0.0138888888888889	0.000160751028806584\\
0.0104166666666667	6.78168402777778e-05\\
0.00757575757575758	2.60873194757492e-05\\
0.00555555555555556	1.02880658436214e-05\\
0.00409836065573771	4.13030165520462e-06\\
0.00297619047619048	1.5817338840298e-06\\
};
\addlegendentry{$h^3$}

\addplot [color=red,line width=1pt,solid]
  table[row sep=crcr]{%
0.05	0.00125\\
0.0357142857142857	0.000325385256143274\\
0.025	7.8125e-05\\
0.0192307692307692	2.73537341129512e-05\\
0.0138888888888889	7.44217725956409e-06\\
0.0104166666666667	2.35475139853395e-06\\
0.00757575757575758	6.58770693832051e-07\\
0.00555555555555556	1.90519737844841e-07\\
0.00409836065573771	5.64248860000631e-08\\
0.00297619047619048	1.56918044050576e-08\\
};
\addlegendentry{$h^4$}

\end{axis}

\end{tikzpicture}%
  }
  \subfigure[first-order jump conditions]{ \label{fig:firstOrder}
		\setlength\figureheight{0.33\linewidth} 
		\setlength\figurewidth{0.33\linewidth} 
		\tikzset{external/export next=false}
\begin{tikzpicture}

\begin{axis}[%
width=0.951\figurewidth,
height=\figureheight,
at={(0\figurewidth,0\figureheight)},
scale only axis,
xmode=log,
xmin=0.00297619047619048,
xmax=0.1,
xminorticks=true,
xlabel style={font=\color{white!15!black}},
xlabel={\scriptsize$h$},
ymode=log,
ymin=1e-06,
ymax=10,
yminorticks=true,
ytick = {1e-8,1e-7,1e-6,1e-5, 1e-4 ,1e-3, 1e-2, 1e-1,1},
ylabel={\scriptsize$E_1$},
axis background/.style={fill=white},
legend style={at={(0.615,0.5)},anchor=north west,legend cell align=left,align=left,draw=white!15!black,draw=none,fill=none},
legend style={font=\scriptsize},
ylabel style={yshift=-5pt},xlabel style={yshift=2.5pt},tick label style={font=\tiny} 
]
\addplot [color=black,line width=1pt,solid,mark=o,mark options={solid}]
  table[row sep=crcr]{%
0.05	0.509869275318805\\
0.0357142857142857	0.293788560525675\\
0.025	0.129948409995766\\
0.0192307692307692	0.0748675186078586\\
0.0138888888888889	0.0367793603540108\\
0.0104166666666667	0.0181048778042401\\
0.00757575757575758	0.00754319120384234\\
0.00555555555555556	0.00325345053004169\\
0.00409836065573771	0.00134333472502977\\
0.00297619047619048	0.000510899652698514\\
};
\addlegendentry{Yee}

\addplot [color=blue,line width=1pt,solid,mark=square,mark options={solid}]
  table[row sep=crcr]{%
0.05	0.239139419629925\\
0.0357142857142857	0.0856530374074441\\
0.025	0.0295588891046585\\
0.0192307692307692	0.011997797777369\\
0.0138888888888889	0.00382484450504707\\
0.0104166666666667	0.00143395741005771\\
0.00757575757575758	0.000472811720963989\\
0.00555555555555556	0.000165254593075241\\
0.00409836065573771	5.4723842986085e-05\\
0.00297619047619048	1.90479449947044e-05\\
};
\addlegendentry{4$^{th}$}

\addplot [color=red,line width=1pt,densely dashed]
  table[row sep=crcr]{%
0.05	6.625\\
0.0357142857142857	2.41435860058309\\
0.025	0.828125\\
0.0192307692307692	0.376934456076468\\
0.0138888888888889	0.141996742112483\\
0.0104166666666667	0.0599048755787037\\
0.00757575757575758	0.0230437988702452\\
0.00555555555555556	0.0090877914951989\\
0.00409836065573771	0.00364843312876408\\
0.00297619047619048	0.00139719826422633\\
};
\addlegendentry{$h^3$}

\addplot [color=red,line width=1pt,solid]
  table[row sep=crcr]{%
0.05	0.111803398874989\\
0.0357142857142857	0.0344355387216861\\
0.025	0.00988211768802619\\
0.0192307692307692	0.00394501163678617\\
0.0138888888888889	0.00126297936168722\\
0.0104166666666667	0.00046143515180106\\
0.00757575757575758	0.000151373980816862\\
0.00555555555555556	5.11218101851804e-05\\
0.00409836065573771	1.76276979038025e-05\\
0.00297619047619048	5.75271051812942e-06\\
};
\addlegendentry{$h^{3.5}$}

\end{axis}
\end{tikzpicture}%
		}		
  \caption{Convergence plots of jump conditions for a scattering of a magnetic dielectric cylinder problem with $\mu^+=1$ and $\mu^-=2$ using the proposed CFM-FDTD schemes.}
   \label{fig:convPlotMagneticJumpCdns}
\end{figure}
Convergence plots of zeroth and first order jump conditions at $t_f = 1$ are shown in 
	Figure~\ref{fig:convPlotMagneticJumpCdns} for both schemes.
A third-order convergence is observed for zeroth and first order jump conditions when the CFM-Yee scheme is used. 
As for the CFM-$4^{th}$ scheme, 
	a fourth-order convergence is obtained for zeroth-order jump conditions while 
	a three and a half order convergence is observed for first-order jump conditions.
%We observe a third and fourth order convergence of zeroth-order jump conditions 
%	for respectively the CFM-Yee scheme and the CFM-$4^{th}$ scheme as expected. 
%A third, 
%	and three and a half order convergence of first-order jump conditions are also observed for both schemes.
According to numerical results, 
	approximations of correction functions coming from minimization problem \eqref{eq:minPblm}
	are consistent with high-order explicit jump conditions coming from Maxwell's equation 
	\eqref{eq:pblmDefinition} and therefore are appropriate to correct FD approximations in the vicinity of the interface. 

%%%%%%%%%%%%%%%%%%%%%%%%%%%%%%%%%%%%%%%%%%%%%%
\subsection{Problems with a Manufactured Solution}
To our knowledge, 
	there is no analytic solution for Maxwell's interface problems with an arbitrary geometry 
	of the interface. 
In order to verify the proposed numerical strategy,
%	we construct manufactured solutions.
%Hence,
	general interface conditions,
	given by 
\begin{subequations} \label{eq:generalInterfCdns}
    \begin{align}
    	\hat{\mathbold{n}}\times\llbracket \mathbold{E} \rrbracket =&\,\, \mathbold{a}(\mathbold{x},t) \quad \text{on } \Gamma \times I ,\label{eq:generalTangentEInterf}\\
    	\hat{\mathbold{n}}\times\llbracket \mathbold{H} \rrbracket =&\,\, \mathbold{b}(\mathbold{x},t) \quad \text{on } \Gamma \times I ,\label{eq:generalTangentHInterf}\\
    	\hat{\mathbold{n}}\cdot\llbracket \epsilon(\mathbold{x})\,\mathbold{E} \rrbracket =&\,\, c(\mathbold{x},t) \quad \text{on } \Gamma \times I ,\label{eq:generalNormalEInterf}\\
    	\hat{\mathbold{n}}\cdot\llbracket \mu(\mathbold{x})\,\mathbold{H} \rrbracket =&\,\, d(\mathbold{x},t) \quad \text{on } \Gamma \times I ,\label{eq:generalNormalHInterf}
    \end{align}
\end{subequations}	
	are considered.
Hence, 
	both tangential and normal components of electromagnetic fields
	can be discontinuous on the interface. 
Moreover, 
	electromagnetic fields are at divergence-free in each subdomain, 
	but not necessarily in the whole domain because of interface conditions 	
	\eqref{eq:generalNormalEInterf} and \eqref{eq:generalNormalHInterf}.
Source terms in each subdomain are given by 
	$\mathbold{f}^+_1(\mathbold{x},t)$ and $\mathbold{f}^-_1(\mathbold{x},t)$ for 
	Faraday's law \eqref{eq:Faraday},
	and by $\mathbold{f}^+_2(\mathbold{x},t)$ and $\mathbold{f}^-_2(\mathbold{x},t)$ for 
	Amp\`{e}re-Maxwell's law \eqref{eq:AmpereMaxwell}.  
It is worth mentioning that these source terms and interface conditions are not substantiated by 	
	physics.
Nevertheless, 
	they can be used to construct manufactured solutions that are needed to verify the 
	proposed numerical strategy for arbitrary complex interfaces.
	
\begin{figure} 
 \centering
  \subfigure[circular]{\label{fig:circleInterfaceGeo}
		\setlength\figureheight{0.21\linewidth} 
		\setlength\figurewidth{0.21\linewidth} 
		\tikzset{external/export next=false}
\begin{tikzpicture}

\begin{axis}[%
width=\figurewidth,
height=\figureheight,
at={(0\figurewidth,0\figureheight)},
scale only axis,
xmin=0,
xmax=1,
xminorticks=true,
xtick = {0, 0.5 ,1},
ymin=0,
ymax=1,
yminorticks=true,
ytick = {0, 0.5 ,1},
axis background/.style={fill=white},
legend style={at={(0.01,0.99)},anchor=north west,legend cell align=left,align=left,draw=white!15!black,draw=none,fill=none},
legend style={font=\scriptsize},
ylabel style={yshift=-5pt},xlabel style={yshift=2.5pt},tick label style={font=\tiny} 
]
\addplot [color=black,line width=1pt,solid]
  table[row sep=crcr]{%
0.75	0.5\\
0.749496669117971	0.515855979914141\\
0.747988703207699	0.531648113393437\\
0.745482174315677	0.547312811090103\\
0.741987175349089	0.56278699679527\\
0.737517779435236	0.578008361424622\\
0.732091983254018	0.592915613915082\\
0.725731634571655	0.607448728022293\\
0.718462344267446	0.621549184025117\\
0.710313383207795	0.635160204363899\\
0.701317564382765	0.64822698226366\\
0.691511110779744	0.660696902421635\\
0.680933509526268	0.672519752870528\\
0.669627352889283	0.683647927164383\\
0.657638166771131	0.694036616072939\\
0.6450142273928	0.703643988012584\\
0.631806366902626	0.712431357487379\\
0.618067768693171	0.720363340861895\\
0.603853753250472	0.72740799883863\\
0.589221555397968	0.733536965066277\\
0.574230093832069	0.738725560361018\\
0.558939733877357	0.742952892080885\\
0.543412044416733	0.746201938253052\\
0.527709549975253	0.748459616115314\\
0.511895478955936	0.749716834795752\\
0.496033509041298	0.749968531918469\\
0.480187510785803	0.749213693987986\\
0.464421290431679	0.747455360470233\\
0.448798332983702	0.744700611553695\\
0.433381546577491	0.740960539639986\\
0.418233009170645	0.736250204678667\\
0.403413718576718	0.730588573526145\\
0.388983346848556	0.723998443572834\\
0.375	0.71650635094611\\
0.361519984033472	0.708142463658693\\
0.348597578215583	0.698940460132708\\
0.336284816513679	0.688937393588565\\
0.32463127807342	0.678173542844716\\
0.313683887581061	0.666692250129073\\
0.303486726314303	0.654539746555151\\
0.294080854642542	0.641764965965693\\
0.285504146691256	0.628419347893352\\
0.277791137836269	0.614556630431853\\
0.270972885641983	0.600232633851653\\
0.265076844803523	0.585505035831417\\
0.260126756596376	0.570433139210357\\
0.256142553278648	0.555077633196635\\
0.253140277830901	0.539500348993337\\
0.251132019356729	0.523764010826046\\
0.250125864404204	0.507931983374517\\
0.250125864404204	0.492068016625483\\
0.251132019356729	0.476235989173954\\
0.253140277830901	0.460499651006662\\
0.256142553278648	0.444922366803365\\
0.260126756596376	0.429566860789643\\
0.265076844803523	0.414494964168583\\
0.270972885641983	0.399767366148347\\
0.277791137836269	0.385443369568148\\
0.285504146691256	0.371580652106648\\
0.294080854642542	0.358235034034307\\
0.303486726314303	0.345460253444849\\
0.313683887581061	0.333307749870927\\
0.32463127807342	0.321826457155284\\
0.336284816513679	0.311062606411435\\
0.348597578215583	0.301059539867292\\
0.361519984033472	0.291857536341307\\
0.375	0.28349364905389\\
0.388983346848557	0.276001556427166\\
0.403413718576718	0.269411426473855\\
0.418233009170645	0.263749795321333\\
0.433381546577491	0.259039460360015\\
0.448798332983702	0.255299388446305\\
0.464421290431679	0.252544639529767\\
0.480187510785803	0.250786306012014\\
0.496033509041298	0.250031468081531\\
0.511895478955936	0.250283165204248\\
0.527709549975253	0.251540383884686\\
0.543412044416733	0.253798061746948\\
0.558939733877357	0.257047107919115\\
0.574230093832069	0.261274439638982\\
0.589221555397968	0.266463034933723\\
0.603853753250472	0.27259200116137\\
0.61806776869317	0.279636659138104\\
0.631806366902626	0.287568642512621\\
0.645014227392799	0.296356011987416\\
0.657638166771131	0.305963383927061\\
0.669627352889283	0.316352072835617\\
0.680933509526267	0.327480247129472\\
0.691511110779744	0.339303097578365\\
0.701317564382765	0.35177301773634\\
0.710313383207795	0.364839795636101\\
0.718462344267446	0.378450815974883\\
0.725731634571655	0.392551271977707\\
0.732091983254018	0.407084386084918\\
0.737517779435236	0.421991638575378\\
0.741987175349089	0.43721300320473\\
0.745482174315677	0.452687188909897\\
0.747988703207699	0.468351886606563\\
0.749496669117971	0.484144020085859\\
0.75	0.5\\
};
\end{axis}
\end{tikzpicture}%
		} 
  \subfigure[5-star]{\label{fig:starInterfaceGeo}
		\setlength\figureheight{0.21\linewidth} 
		\setlength\figurewidth{0.21\linewidth} 
		\tikzset{external/export next=false}
\begin{tikzpicture}

\begin{axis}[%
width=\figurewidth,
height=\figureheight,
at={(0\figurewidth,0\figureheight)},
scale only axis,
xmin=0,
xmax=1,
xminorticks=true,
xtick = {0, 0.5 ,1},
ymin=0,
ymax=1,
yminorticks=true,
ytick = {0, 0.5 ,1},
axis background/.style={fill=white},
legend style={at={(0.01,0.99)},anchor=north west,legend cell align=left,align=left,draw=white!15!black,draw=none,fill=none},
legend style={font=\scriptsize},
ylabel style={yshift=-5pt},xlabel style={yshift=2.5pt},tick label style={font=\tiny} 
]
\addplot [color=black,line width=1pt,solid]
  table[row sep=crcr]{%
0.75	0.5\\
0.765066930189006	0.516845499123648\\
0.777395596897263	0.535400996867555\\
0.785474949486592	0.555020786717682\\
0.788201994577402	0.574778085589939\\
0.785015355939205	0.593608069890764\\
0.77595740609245	0.610476680159308\\
0.76165735877021	0.624549447537036\\
0.743238589707334	0.6353343167743\\
0.722163808645616	0.64277601035563\\
0.700040086323819	0.647286395096484\\
0.678410979268275	0.649704586901987\\
0.658564374561849	0.651190825763033\\
0.641382028781861	0.653067981559738\\
0.627250504584707	0.656632481897858\\
0.616044232300894	0.662961322323235\\
0.607181039780345	0.672742897877369\\
0.599740192449566	0.686156495288\\
0.592624237639838	0.702818789561018\\
0.584739992587695	0.721806497324522\\
0.57517168653353	0.74175374251308\\
0.563320871123122	0.761012185776997\\
0.54899298926918	0.777853049331026\\
0.532418651828539	0.790684106941049\\
0.514207511767961	0.798252376567682\\
0.495242705950108	0.799804996894593\\
0.47653268388508	0.795186358147882\\
0.459043571509212	0.784858217139645\\
0.443538105232274	0.769840045925559\\
0.430446197184106	0.751577768613709\\
0.41978749849487	0.731758802741356\\
0.411158596481512	0.712098573477648\\
0.403787839421047	0.694127399905181\\
0.396650635094611	0.67900635094611\\
0.388628899709452	0.667396393140633\\
0.378691488291971	0.659397523849207\\
0.366068930318577	0.654564693923155\\
0.350396159890447	0.651996581389543\\
0.331801144729265	0.650483204539873\\
0.310924802629974	0.64869037868741\\
0.28886729267282	0.645354240940172\\
0.267066301407991	0.639458144361026\\
0.247122803447367	0.63036725568928\\
0.230597541587669	0.617902712298803\\
0.218806015883607	0.602346240273091\\
0.212640420665489	0.584376385492237\\
0.212443602463266	0.564947476486013\\
0.217953180803277	0.545130682692334\\
0.228324437528831	0.525941870830668\\
0.24222977195585	0.508182636263719\\
0.258021956852558	0.492318669514685\\
0.273939601184627	0.478413849178577\\
0.288327374858526	0.466129984705659\\
0.299841504094031	0.454792210092742\\
0.307613092527262	0.443510107071522\\
0.311347673723439	0.431336168610256\\
0.311348229696296	0.417437444595497\\
0.308459472225172	0.401253994825574\\
0.303941991974521	0.382619448574323\\
0.299294416612264	0.361824309008787\\
0.296048649998632	0.339610885577108\\
0.295566630432858	0.317098704281728\\
0.298866396256393	0.295649495700111\\
0.30650070270878	0.276689906746026\\
0.318503668139195	0.261516603583791\\
0.334411068357493	0.251111465823247\\
0.353349364905389	0.24599364905389\\
0.374178854276066	0.246130512759513\\
0.395668840671924	0.250921426425358\\
0.416678519846419	0.259258393384021\\
0.436316895970877	0.269656689333738\\
0.454058560735131	0.280438822818169\\
0.469799009354145	0.289947496199178\\
0.483842337686526	0.296758970171911\\
0.496824312132489	0.299867933057655\\
0.50958344614391	0.298818706976178\\
0.523000448121967	0.293764874710422\\
0.537831099564285	0.285449172824922\\
0.554558596631592	0.275106401615227\\
0.573288501130608	0.264302621791043\\
0.593703118208241	0.254732567191969\\
0.615083268861105	0.248002791883759\\
0.636395344936775	0.245429813564209\\
0.656431694024907	0.247880182902612\\
0.673984222484705	0.255673346298067\\
0.688025828957554	0.268559249751979\\
0.697872676996705	0.285772127230971\\
0.703302644490686	0.306151320021977\\
0.704611242291214	0.328310782058717\\
0.702595042441711	0.350832430569164\\
0.698462957769975	0.372455601627831\\
0.693686098827558	0.392235948724066\\
0.689805910373101	0.40965199149245\\
0.688226560415586	0.424645452329144\\
0.690020202931268	0.437591347041521\\
0.695772356120776	0.449204091999399\\
0.705489399144761	0.460395164537477\\
0.718581809518134	0.47210477008068\\
0.733926408046936	0.485133539295366\\
0.75	0.5\\
};
\end{axis}
\end{tikzpicture}%
		} 
  \subfigure[3-star]{\label{fig:triStarInterfaceGeo}
		\setlength\figureheight{0.21\linewidth} 
		\setlength\figurewidth{0.21\linewidth} 
		\tikzset{external/export next=false}
\begin{tikzpicture}

\begin{axis}[%
width=\figurewidth,
height=\figureheight,
at={(0\figurewidth,0\figureheight)},
scale only axis,
xmin=0,
xmax=1,
xminorticks=true,
xtick = {0, 0.5 ,1},
ymin=0,
ymax=1,
yminorticks=true,
ytick = {0, 0.5 ,1},
axis background/.style={fill=white},
legend style={at={(0.01,0.99)},anchor=north west,legend cell align=left,align=left,draw=white!15!black,draw=none,fill=none},
legend style={font=\scriptsize},
ylabel style={yshift=-5pt},xlabel style={yshift=2.5pt},tick label style={font=\tiny} 
]
\addplot [color=black,line width=1pt,solid]
  table[row sep=crcr]{%
0.75	0.55\\
0.768383691155712	0.567056285485863\\
0.784855939371776	0.586353079609766\\
0.798569247672799	0.607544505835653\\
0.808783283651588	0.630118192268234\\
0.814908287960172	0.653425855534761\\
0.816539300787274	0.676723220029148\\
0.813479080096982	0.699216694779779\\
0.805748304479489	0.720113788035931\\
0.793582461693234	0.738673991711258\\
0.777415669648107	0.754256830101566\\
0.757852505596638	0.766363942344277\\
0.735629678650627	0.774672444790369\\
0.711570021888213	0.779057373755364\\
0.686531762311539	0.779601705541283\\
0.661356319017503	0.776593244583852\\
0.636817963590651	0.770508511213973\\
0.61357854671968	0.761984593953246\\
0.592150159874138	0.75178070405503\\
0.57286808220072	0.740731832587438\\
0.555875694011952	0.729697420209933\\
0.541122262361211	0.719508274169579\\
0.528373671098689	0.710915085057808\\
0.517235331010121	0.704541799856712\\
0.507185718903124	0.700846803803962\\
0.497618308349223	0.700094370975715\\
0.48788911327206	0.702338193667119\\
0.477366703462387	0.707418035147949\\
0.465481391646248	0.714969718104235\\
0.4517703451712	0.724447815239603\\
0.435915638275791	0.735159603166056\\
0.417772716391469	0.74630812733257\\
0.397387362750012	0.757041649900452\\
0.375	0.76650635094611\\
0.351036977897187	0.773898951759353\\
0.326089339866257	0.778515940109955\\
0.300880372261704	0.779796300355625\\
0.276223968364629	0.777355071772135\\
0.252976437718672	0.771005649468394\\
0.231984821811756	0.760769465763369\\
0.214035011603472	0.746872499529074\\
0.199802985787138	0.729728904825291\\
0.189812299298019	0.709912874076284\\
0.184400561411509	0.68812060227235\\
0.183697076668585	0.66512484910402\\
0.18761311595587	0.641725065201822\\
0.195845484459602	0.61869632675674\\
0.207893209093245	0.596740392003832\\
0.223086334355692	0.576442049038377\\
0.240625045742688	0.558233576556572\\
0.259626683065719	0.542369609807538\\
0.279177704357766	0.528914027386286\\
0.298387346568558	0.517739694017157\\
0.316439622097695	0.50854106036347\\
0.332640397236881	0.500858786781107\\
0.34645661293846	0.494114777441185\\
0.357545209872456	0.487655334569043\\
0.365769976374519	0.480799613212578\\
0.371205307595373	0.472890209038588\\
0.374126697681612	0.463342567597689\\
0.374988630816851	0.451689972653066\\
0.374391337443451	0.437621149210248\\
0.37303858778221	0.421007986082703\\
0.371689260765654	0.401921513178496\\
0.37110581656491	0.380635019844539\\
0.372002990169758	0.357614024441967\\
0.375	0.33349364905389\\
0.380579330947101	0.309044762754784\\
0.389054720761967	0.285130980280279\\
0.400550380065498	0.262659193808722\\
0.414992747983782	0.242526735959632\\
0.432115274321156	0.225568494996845\\
0.45147587740097	0.212507314207483\\
0.472485908299546	0.203910805691148\\
0.494448709733373	0.200157307138778\\
0.516605239008747	0.201413134212458\\
0.538183768940384	0.207622567626085\\
0.558450417734776	0.218511208551704\\
0.576757205393503	0.233602490007808\\
0.592584493652186	0.252246299487896\\
0.605575028595216	0.273657902454885\\
0.615557346626805	0.29696470637777\\
0.622556990666661	0.321257912229455\\
0.6267947702146	0.345645796239216\\
0.628672135768096	0.369305268558684\\
0.628744571230722	0.391528473395405\\
0.627684683890353	0.411761519426597\\
0.626237340401908	0.429632939049313\\
0.625169715962851	0.444970137501007\\
0.625219459117422	0.457802865574245\\
0.627044304722356	0.468353582983459\\
0.631176384055404	0.477015419985697\\
0.637984189046328	0.484319238735193\\
0.647644665720762	0.490891992198984\\
0.660127270910301	0.497409132685517\\
0.67519106704659	0.504544198677694\\
0.692395100958555	0.512918883655448\\
0.711121467043621	0.523056852822891\\
0.730609647080229	0.535344325657581\\
0.75	0.55\\
};
\end{axis}
\end{tikzpicture}%
		} 
  \caption{Different geometries of the interface.}
  \label{fig:interfaceGeo}
\end{figure}
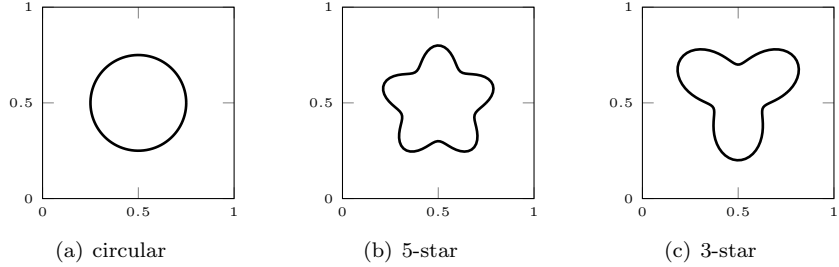

The domain is $\Omega = [0,1]\times[0,1]$ and the time interval is $I = [0,1]$. 
The physical parameters are given by $\mu^+ = 2$, 
	$\epsilon^+ =1$, 
	$\mu^- = \sin(5\,\pi\,x\,y)+2$ and $\epsilon^- = 2\,e^{x\,y}$.
The magnetic permeability and the electrical permittivity have been chosen
	in such a way that electromagnetic fields are at
	divergence-free in each subdomain.
The manufactured solutions are :
\begin{equation*}
	\begin{aligned}
		H_x^+ =&\,\, 0.5\,\sin(2\,\pi\,x)\,\sin(2\,\pi\,y)\,\sin(2\,\pi\,t), \\
		H_y^+ =&\,\, 0.5\,\cos(2\,\pi\,x)\,\cos(2\,\pi\,y)\,\sin(2\,\pi\,t), \\
		E_z^+ =&\,\, \sin(2\,\pi\,x)\,\cos(2\,\pi\,y)\,\cos(2\,\pi\,t)
	\end{aligned}
\end{equation*}
	in $\Omega^+$, 
	and 
\begin{equation*}
	\begin{aligned}
		H_x^- =&\,\, -x\,e^{-x\,y}\,\sin(2\,\pi\,t), \\
		H_y^- =&\,\, y\,e^{-x\,y}\,\sin(2\,\pi\,t), \\
		E_z^- =&\,\, \sin(2\,\pi\,x\,y)\,\cos(2\,\pi\,t)
	\end{aligned}
\end{equation*}
	in $\Omega^-$.
The associated source terms are $\mathbold{f}_1^+ = 0$, $f_2^+=0$
	and 
\begin{equation*}
	\begin{aligned}
		f_{1,x}^- =&\,\, 2\,\pi\,x\,\big(\cos(2\,\pi\,x\,y) - (\sin(5\,\pi\,x\,y)+2)\,e^{-x\,y}\big)\,\cos(2\,\pi\,t), \\
		f_{1,y}^- =&\,\, 2\,\pi\,y\,\big((\sin(5\,\pi\,x\,y)+2)\,e^{-x\,y}-\cos(2\,\pi\,x\,y)\big)\,\cos(2\,\pi\,t), \\
		f_2^-  =&\,\, \big((x^2+y^2)\,e^{-x\,y}-4\,\pi\,e^{x\,y}\,\sin(2\,\pi\,x\,y) \big)\,\sin(2\,\pi\,t).
	\end{aligned}
\end{equation*}
We consider geometries of the interface that are illustrated in Figure~\ref{fig:interfaceGeo}.
%**********************************************
%**********************************************
Periodic boundary conditions are imposed on all $\partial \Omega$ for both CFM-FDTD schemes. 
The mesh grid size is $h=\Delta x= \Delta y$  and the time step is $\Delta t = \tfrac{h}{2}$ with
	$h \in \big\{\tfrac{1}{20},\tfrac{1}{28},\tfrac{1}{40},\tfrac{1}{52}, \tfrac{1}{72}, \tfrac{1}{96}, \tfrac{1}{132}, \tfrac{1}{180},\tfrac{1}{244},\tfrac{1}{336},\tfrac{1}{460}\big\}$.
For local patches, 
	we choose $\ell_h = \beta\,h$ with $\beta=8$ for the 5-star interface and $\beta = 7$ for either
	the circular or 3-star interface.	
All other parameters are the same as in subsection~\ref{sec:numExpScattering}. 
Figure~\ref{fig:convPlotManuSolPblms} shows convergence plots of $\mathbold{U} = [H_x,H_y,E_z]^T$ for all geometries of the interface.
We observe a second-order convergence in $L^2$-norm for the CFM-Yee scheme.
As for the CFM-$4^{th}$ scheme, 
	the expected order is not clearly observed for smaller mesh grid sizes.
Since the error of $\mathbold{U}$ is already low for this scheme, 
	this suggests a limitation due to the use of double-precision arithmetic and 
	therefore a more accurate floating-point arithmetic should remedy this issue. 
Nevertheless, 
	a global fourth-order convergence is observed using the $L^2$-norm.
These results are in agreement with the theory. 
%**********************************************
 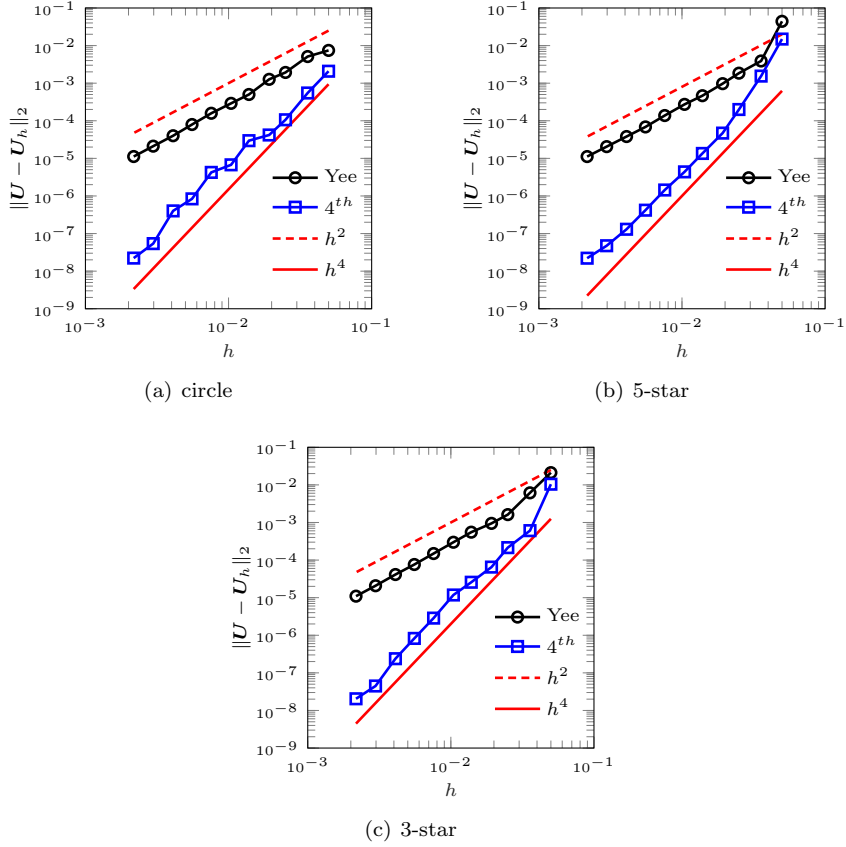
\begin{figure} 
 \centering
  \subfigure[circle]{ \label{fig:convPlotCircularInterface}
		\setlength\figureheight{0.33\linewidth} 
		\setlength\figurewidth{0.33\linewidth} 
		\tikzset{external/export next=false}
\begin{tikzpicture}

\begin{axis}[%
width=0.951\figurewidth,
height=\figureheight,
at={(0\figurewidth,0\figureheight)},
scale only axis,
xmode=log,
xmin=0.001,
xmax=0.1,
xminorticks=true,
xlabel style={font=\color{white!15!black}},
xlabel={\scriptsize$h$},
ymode=log,
ymin=1e-09,
ymax=0.1,
yminorticks=true,
ytick = {1e-9,1e-8,1e-7,1e-6,1e-5, 1e-4 ,1e-3, 1e-2, 1e-1,1},
ylabel={\scriptsize$\|\mathbold{U}-\mathbold{U}_h\|_{2}$},
axis background/.style={fill=white},
legend style={at={(0.62,0.5)},anchor=north west,legend cell align=left,align=left,draw=white!15!black,draw=none,fill=none},
legend style={font=\scriptsize},
ylabel style={yshift=-5pt},xlabel style={yshift=2.5pt},tick label style={font=\tiny} 
]
\addplot [color=black,line width=1pt,solid,mark=o,mark options={solid}]
  table[row sep=crcr]{%
0.05	0.00744290956315461\\
0.0357142857142857	0.00510967592690622\\
0.025	0.00192933742430878\\
0.0192307692307692	0.00126580390863501\\
0.0138888888888889	0.000500405268228166\\
0.0104166666666667	0.000290577949055658\\
0.00757575757575758	0.000158361176612198\\
0.00555555555555556	7.94977756962916e-05\\
0.00409836065573771	3.99924329005771e-05\\
0.00297619047619048	2.11137384499263e-05\\
0.00217391304347826	1.11929635146311e-05\\
};
\addlegendentry{Yee}

\addplot [color=blue,line width=1pt,solid,mark=square,mark options={solid}]
  table[row sep=crcr]{%
0.05	0.00208050810514959\\
0.0357142857142857	0.000548187512932165\\
0.025	0.0001063195801498\\
0.0192307692307692	4.18591710175152e-05\\
0.0138888888888889	2.95771112819686e-05\\
0.0104166666666667	6.73530988244807e-06\\
0.00757575757575758	4.23599787417725e-06\\
0.00555555555555556	8.37329124812881e-07\\
0.00409836065573771	4.00669912397024e-07\\
0.00297619047619048	5.40168248807801e-08\\
0.00217391304347826	2.22486619125347e-08\\
};
\addlegendentry{4$^{th}$}

\addplot [color=red,line width=1pt,densely dashed]
  table[row sep=crcr]{%
0.05	0.025\\
0.0357142857142857	0.0127551020408163\\
0.025	0.00625\\
0.0192307692307692	0.00369822485207101\\
0.0138888888888889	0.00192901234567901\\
0.0104166666666667	0.00108506944444444\\
0.00757575757575758	0.000573921028466483\\
0.00555555555555556	0.000308641975308642\\
0.00409836065573771	0.000167965600644988\\
0.00297619047619048	8.85770975056689e-05\\
0.00217391304347826	4.72589792060492e-05\\
};
\addlegendentry{$h^2$}

\addplot [color=red,line width=1pt,solid]
  table[row sep=crcr]{%
0.05	0.0009375\\
0.0357142857142857	0.000244038942107455\\
0.025	5.859375e-05\\
0.0192307692307692	2.05153005847134e-05\\
0.0138888888888889	5.58163294467307e-06\\
0.0104166666666667	1.76606354890046e-06\\
0.00757575757575758	4.94078020374038e-07\\
0.00555555555555556	1.42889803383631e-07\\
0.00409836065573771	4.23186645000474e-08\\
0.00297619047619048	1.17688533037932e-08\\
0.00217391304347826	3.35011667339668e-09\\
};
\addlegendentry{$h^4$}

\end{axis}

\end{tikzpicture}%
  }
  \subfigure[5-star]{ \label{fig:convPlot5StarInterface}
		\setlength\figureheight{0.33\linewidth} 
		\setlength\figurewidth{0.33\linewidth} 
		\tikzset{external/export next=false}
\begin{tikzpicture}

\begin{axis}[%
width=0.951\figurewidth,
height=\figureheight,
at={(0\figurewidth,0\figureheight)},
scale only axis,
xmode=log,
xmin=0.001,
xmax=0.1,
xminorticks=true,
xlabel style={font=\color{white!15!black}},
xlabel={\scriptsize$h$},
ymode=log,
ymin=1e-09,
ymax=0.1,
yminorticks=true,
ytick = {1e-9,1e-8,1e-7,1e-6,1e-5, 1e-4 ,1e-3, 1e-2, 1e-1,1},
ylabel={\scriptsize$\|\mathbold{U}-\mathbold{U}_h\|_{2}$},
axis background/.style={fill=white},
legend style={at={(0.62,0.5)},anchor=north west,legend cell align=left,align=left,draw=white!15!black,draw=none,fill=none},
legend style={font=\scriptsize},
ylabel style={yshift=-5pt},xlabel style={yshift=2.5pt},tick label style={font=\tiny} 
]
\addplot [color=black,line width=1pt,solid,mark=o,mark options={solid}]
  table[row sep=crcr]{%
0.05	0.0442020684051124\\
0.0357142857142857	0.00389429730217501\\
0.025	0.00184603000824843\\
0.0192307692307692	0.000970990518888939\\
0.0138888888888889	0.00046452556270568\\
0.0104166666666667	0.000272934175733855\\
0.00757575757575758	0.000137830182264932\\
0.00555555555555556	6.8555493277162e-05\\
0.00409836065573771	3.80840926962334e-05\\
0.00297619047619048	2.05212028149254e-05\\
0.00217391304347826	1.11727300783434e-05\\
};
\addlegendentry{Yee}

\addplot [color=blue,line width=1pt,solid,mark=square,mark options={solid}]
  table[row sep=crcr]{%
0.05	0.0149133264838812\\
0.0357142857142857	0.00154349513441347\\
0.025	0.000199375605383135\\
0.0192307692307692	4.69890353494022e-05\\
0.0138888888888889	1.35138599166014e-05\\
0.0104166666666667	4.37370018812271e-06\\
0.00757575757575758	1.43983327067124e-06\\
0.00555555555555556	4.18417808671674e-07\\
0.00409836065573771	1.30083262563218e-07\\
0.00297619047619048	4.74269885332581e-08\\
0.002173913043478	2.211376854674683e-08\\
};
\addlegendentry{4$^{th}$}

\addplot [color=red,line width=1pt,densely dashed]
  table[row sep=crcr]{%
0.05	0.02\\
0.0357142857142857	0.0102040816326531\\
0.025	0.005\\
0.0192307692307692	0.00295857988165681\\
0.0138888888888889	0.00154320987654321\\
0.0104166666666667	0.000868055555555556\\
0.00757575757575758	0.000459136822773186\\
0.00555555555555556	0.000246913580246914\\
0.00409836065573771	0.00013437248051599\\
0.00297619047619048	7.08616780045351e-05\\
0.00217391304347826	3.78071833648393e-05\\
};
\addlegendentry{$h^2$}

\addplot [color=red,line width=1pt,solid]
  table[row sep=crcr]{%
0.05	0.000625\\
0.0357142857142857	0.000162692628071637\\
0.025	3.90625e-05\\
0.0192307692307692	1.36768670564756e-05\\
0.0138888888888889	3.72108862978204e-06\\
0.0104166666666667	1.17737569926698e-06\\
0.00757575757575758	3.29385346916026e-07\\
0.00555555555555556	9.52598689224204e-08\\
0.00409836065573771	2.82124430000316e-08\\
0.00297619047619048	7.84590220252878e-09\\
0.00217391304347826	2.23341111559779e-09\\
};
\addlegendentry{$h^4$}

\end{axis}

\end{tikzpicture}%
		}		
  \subfigure[3-star]{ \label{fig:convPlot3StarInterface}
		\setlength\figureheight{0.33\linewidth} 
		\setlength\figurewidth{0.33\linewidth} 
		\tikzset{external/export next=false}
\begin{tikzpicture}

\begin{axis}[%
width=0.951\figurewidth,
height=\figureheight,
at={(0\figurewidth,0\figureheight)},
scale only axis,
xmode=log,
xmin=0.001,
xmax=0.1,
xminorticks=true,
xlabel style={font=\color{white!15!black}},
xlabel={\scriptsize$h$},
ymode=log,
ymin=1e-09,
ymax=0.1,
yminorticks=true,
ytick = {1e-9,1e-8,1e-7,1e-6,1e-5, 1e-4 ,1e-3, 1e-2, 1e-1,1},
ylabel={\scriptsize$\|\mathbold{U}-\mathbold{U}_h\|_{2}$},
axis background/.style={fill=white},
legend style={at={(0.62,0.5)},anchor=north west,legend cell align=left,align=left,draw=white!15!black,draw=none,fill=none},
legend style={font=\scriptsize},
ylabel style={yshift=-5pt},xlabel style={yshift=2.5pt},tick label style={font=\tiny} 
]
\addplot [color=black,line width=1pt,solid,mark=o,mark options={solid}]
  table[row sep=crcr]{%
0.05	0.0211828748996994\\
0.0357142857142857	0.00614890457284173\\
0.025	0.00162025965486068\\
0.0192307692307692	0.000947023731078167\\
0.0138888888888889	0.000552378004432509\\
0.0104166666666667	0.000298140604136533\\
0.00757575757575758	0.000149310099118916\\
0.00555555555555556	7.62441738105615e-05\\
0.00409836065573771	4.12817128313407e-05\\
0.00297619047619048	2.07745309340761e-05\\
0.00217391304347826	1.09690167377275e-05\\
};
\addlegendentry{Yee}

\addplot [color=blue,line width=1pt,solid,mark=square,mark options={solid}]
  table[row sep=crcr]{%
0.05	0.0104339923066429\\
0.0357142857142857	0.0006070102395971\\
0.025	0.000214386400068605\\
0.0192307692307692	6.5077604733914e-05\\
0.0138888888888889	2.57298425380583e-05\\
0.0104166666666667	1.1784488109891e-05\\
0.00757575757575758	2.85998258185064e-06\\
0.00555555555555556	8.25157615452185e-07\\
0.00409836065573771	2.38117860095706e-07\\
0.00297619047619048	4.47456267937577e-08\\
0.00217391304347826	2.04505822278086e-08\\
};
\addlegendentry{4$^{th}$}

\addplot [color=red,line width=1pt,densely dashed]
  table[row sep=crcr]{%
0.05	0.025\\
0.0357142857142857	0.0127551020408163\\
0.025	0.00625\\
0.0192307692307692	0.00369822485207101\\
0.0138888888888889	0.00192901234567901\\
0.0104166666666667	0.00108506944444444\\
0.00757575757575758	0.000573921028466483\\
0.00555555555555556	0.000308641975308642\\
0.00409836065573771	0.000167965600644988\\
0.00297619047619048	8.85770975056689e-05\\
0.00217391304347826	4.72589792060492e-05\\
};
\addlegendentry{$h^2$}

\addplot [color=red,line width=1pt,solid]
  table[row sep=crcr]{%
0.05	0.00125\\
0.0357142857142857	0.000325385256143274\\
0.025	7.8125e-05\\
0.0192307692307692	2.73537341129512e-05\\
0.0138888888888889	7.44217725956409e-06\\
0.0104166666666667	2.35475139853395e-06\\
0.00757575757575758	6.58770693832051e-07\\
0.00555555555555556	1.90519737844841e-07\\
0.00409836065573771	5.64248860000631e-08\\
0.00297619047619048	1.56918044050576e-08\\
0.00217391304347826	4.46682223119557e-09\\
};
\addlegendentry{$h^4$}

\end{axis}

\end{tikzpicture}%
  }
  \caption{Convergence plots for problems with a manufactured solution using the proposed CFM-FDTD schemes. It is recalled that $\mathbold{U} = [H_x,H_y,E_z]^T$.}
   \label{fig:convPlotManuSolPblms}
\end{figure}
%***************************************
Figure~\ref{fig:plotInterfaceFields} illustrates the computed solutions for different geometries of the 
	interface.
One can observe that there is no spurious oscillation in the vicinity of the interface.
%***************************************
\begin{figure}     
	\centering
	\subfigure[circular]{
	\stackunder[5pt]{
	\stackunder[5pt]{\includegraphics[width=1.58in]{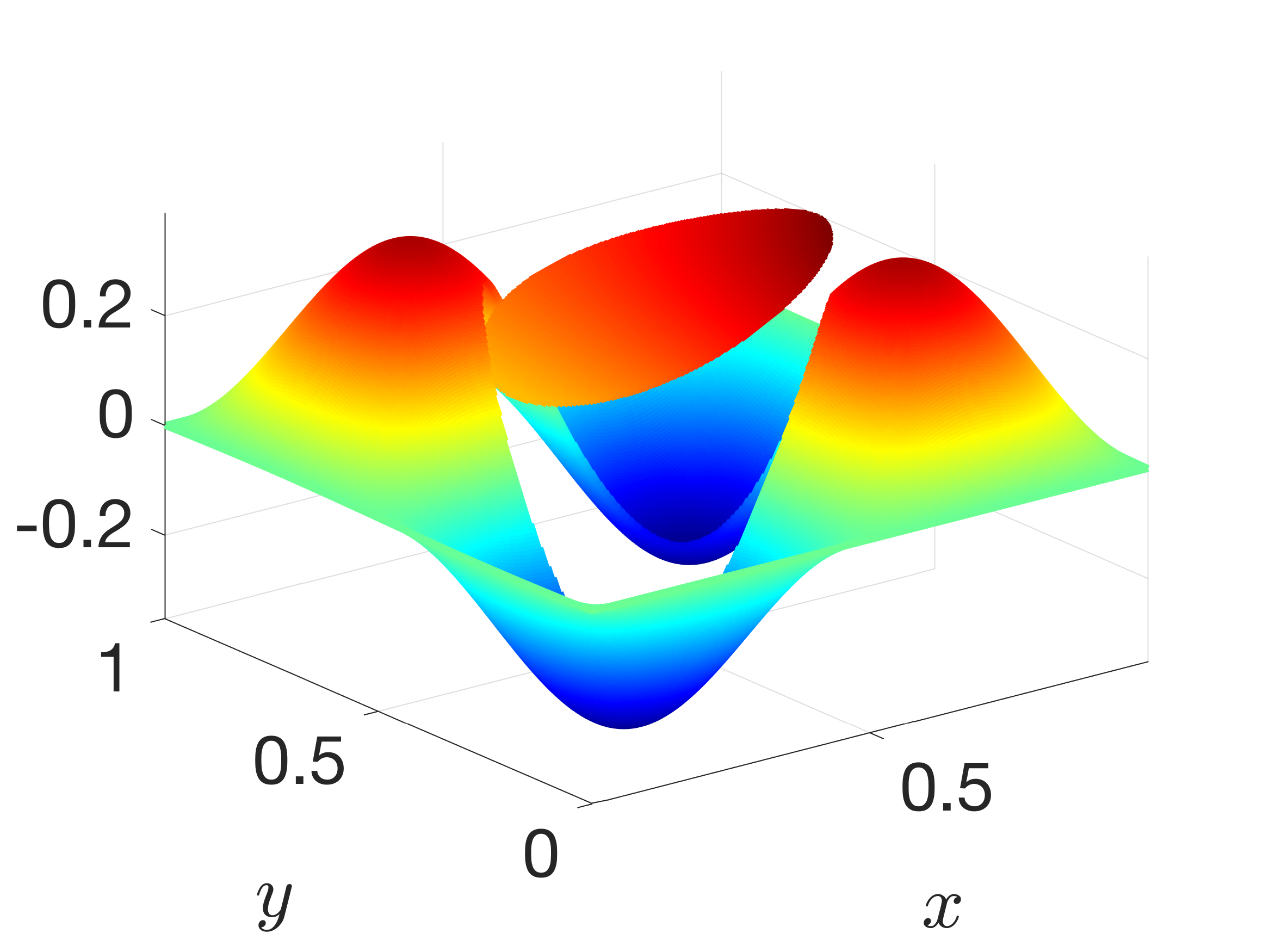}}{\footnotesize$H_x$}
	\stackunder[5pt]{	\includegraphics[width=1.58in]{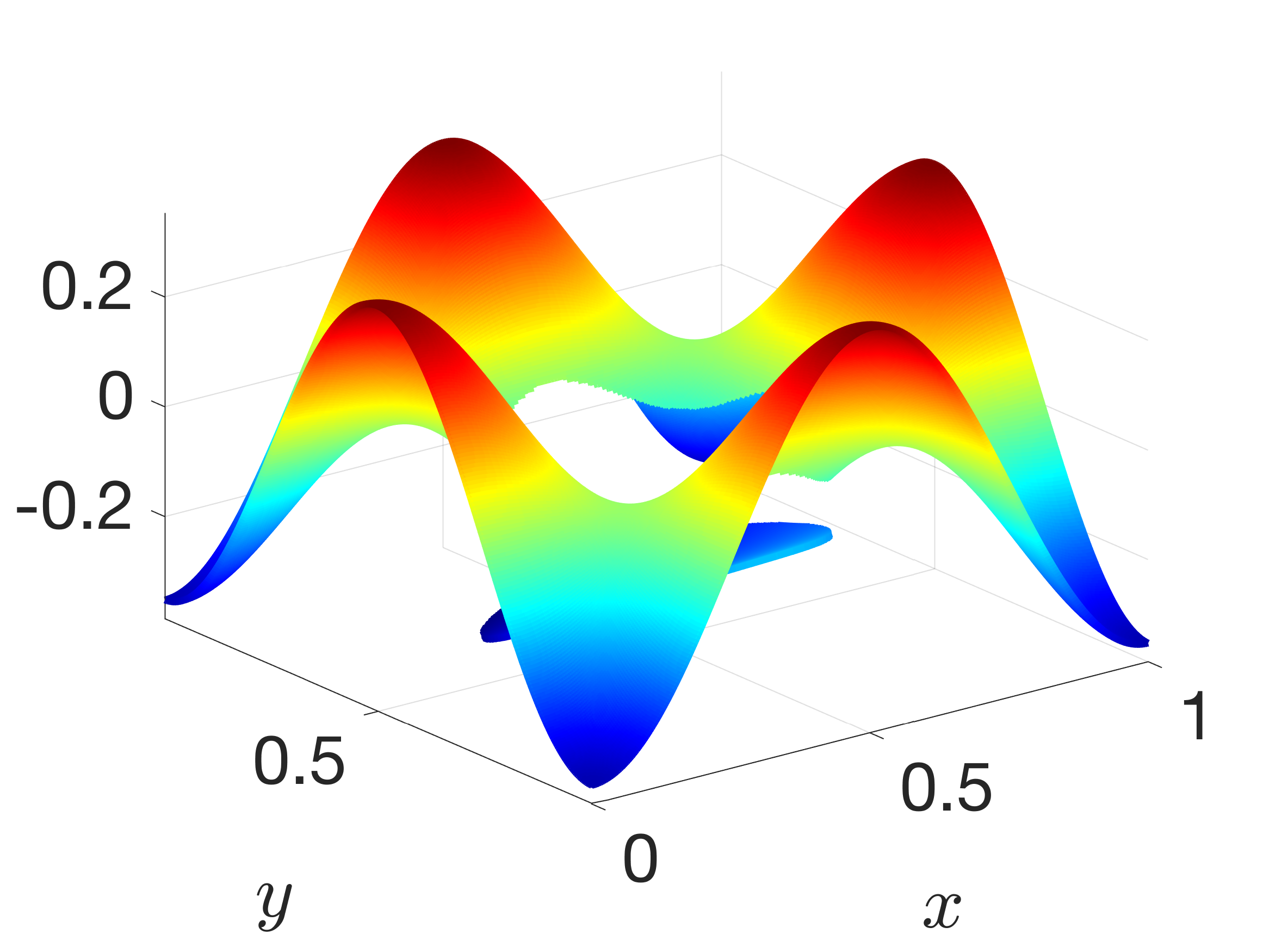}}{\footnotesize$H_y$}
	\stackunder[5pt]{	\includegraphics[width=1.58in]{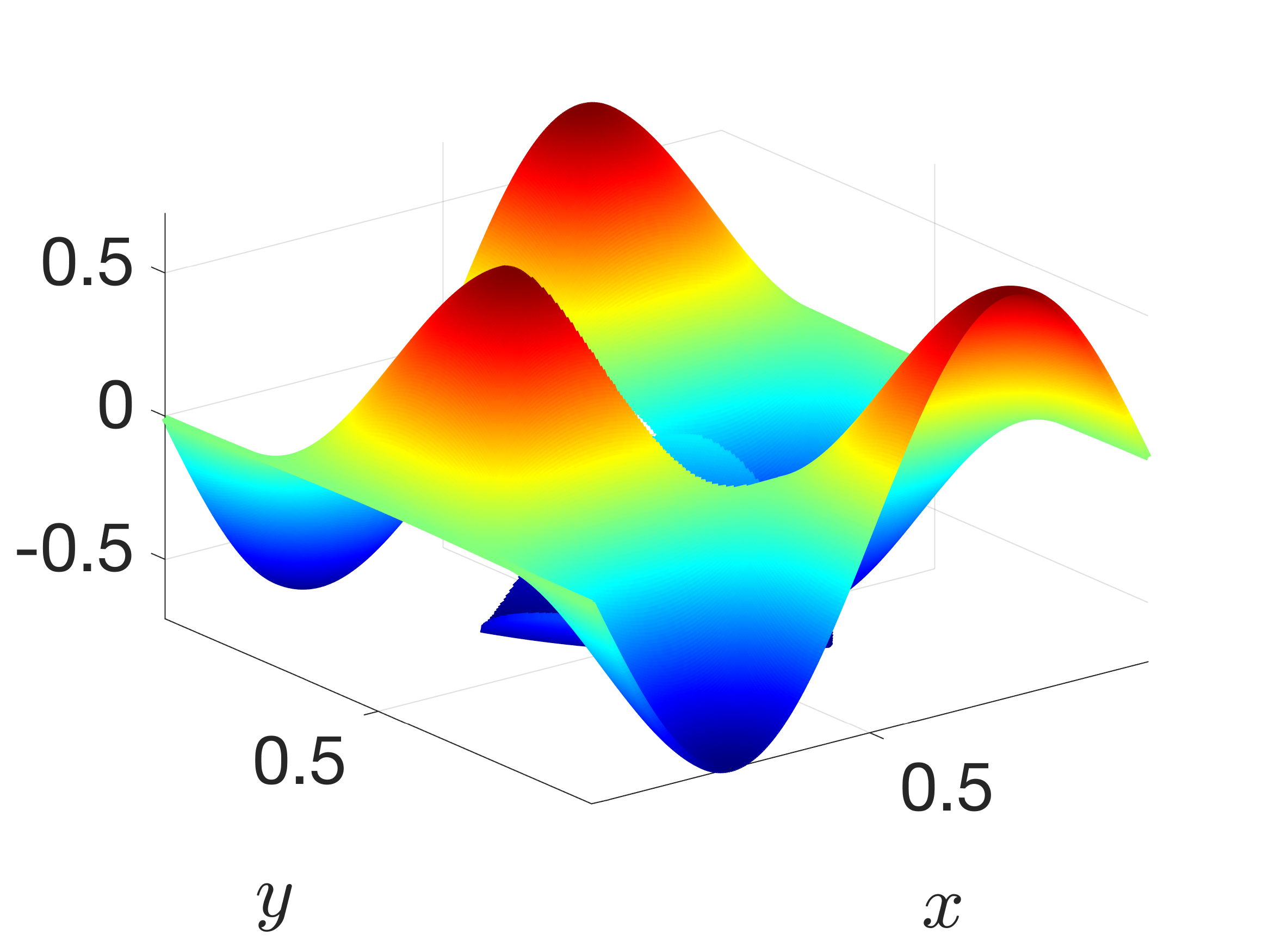}}{\footnotesize$E_z$}
	} 
	{\phantom{A}}
	\label{fig:circleInterfacePEC}
	}
	\subfigure[5-star]{
	\stackunder[5pt]{
	\stackunder[5pt]{\includegraphics[width=1.58in]{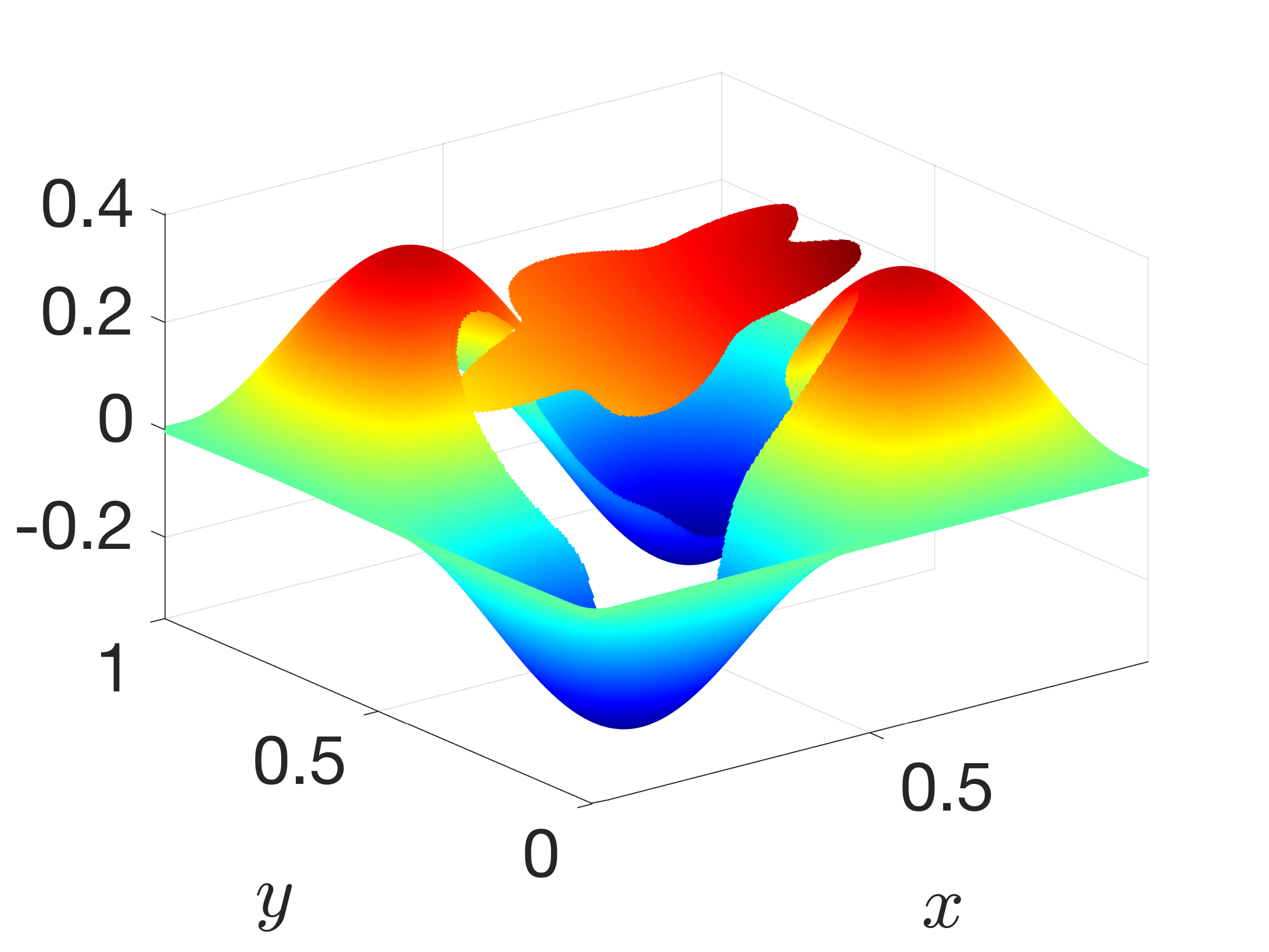}}{\footnotesize$H_x$}
	\stackunder[5pt]{	\includegraphics[width=1.58in]{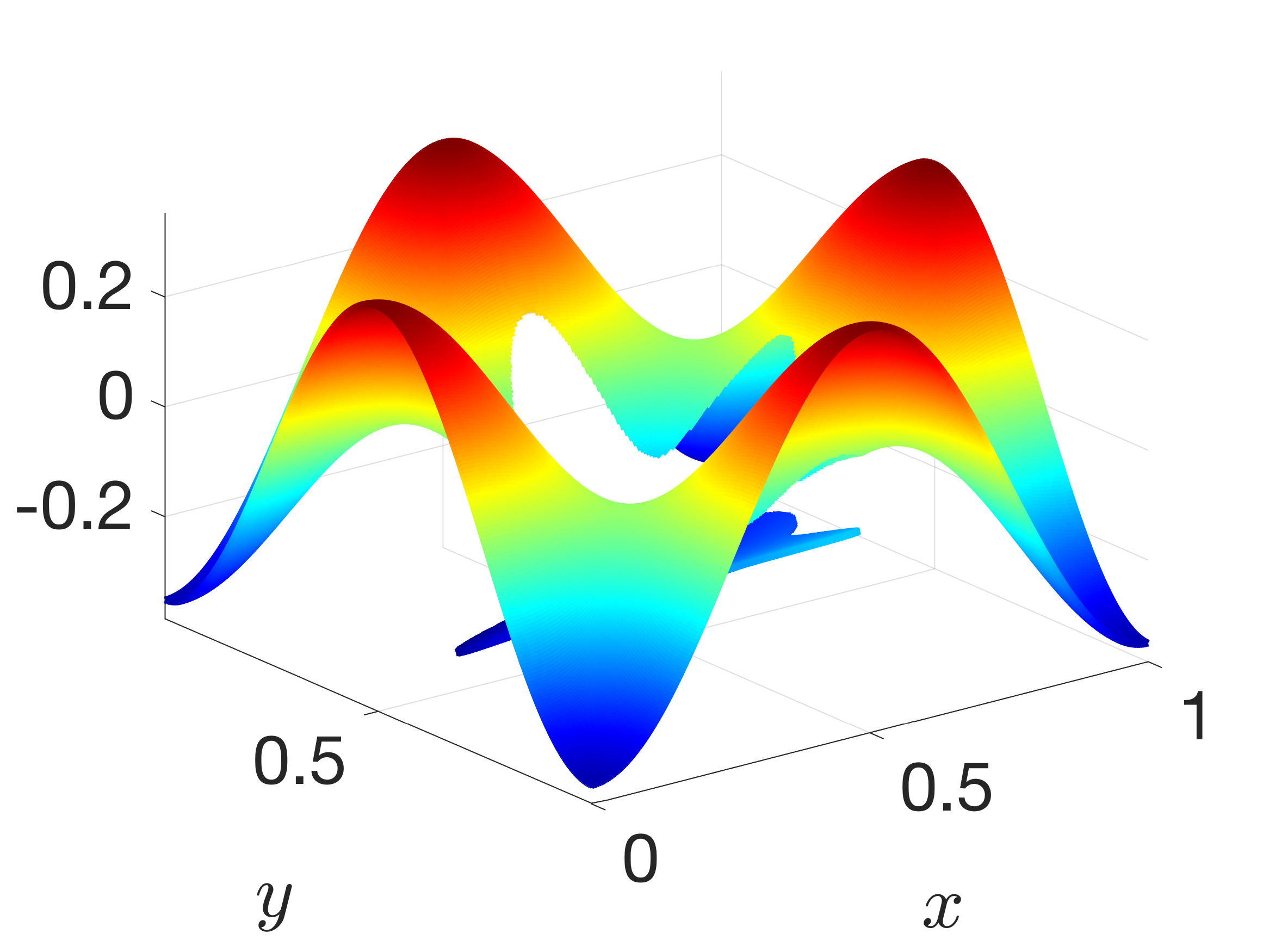}}{\footnotesize$H_y$}
	\stackunder[5pt]{	\includegraphics[width=1.58in]{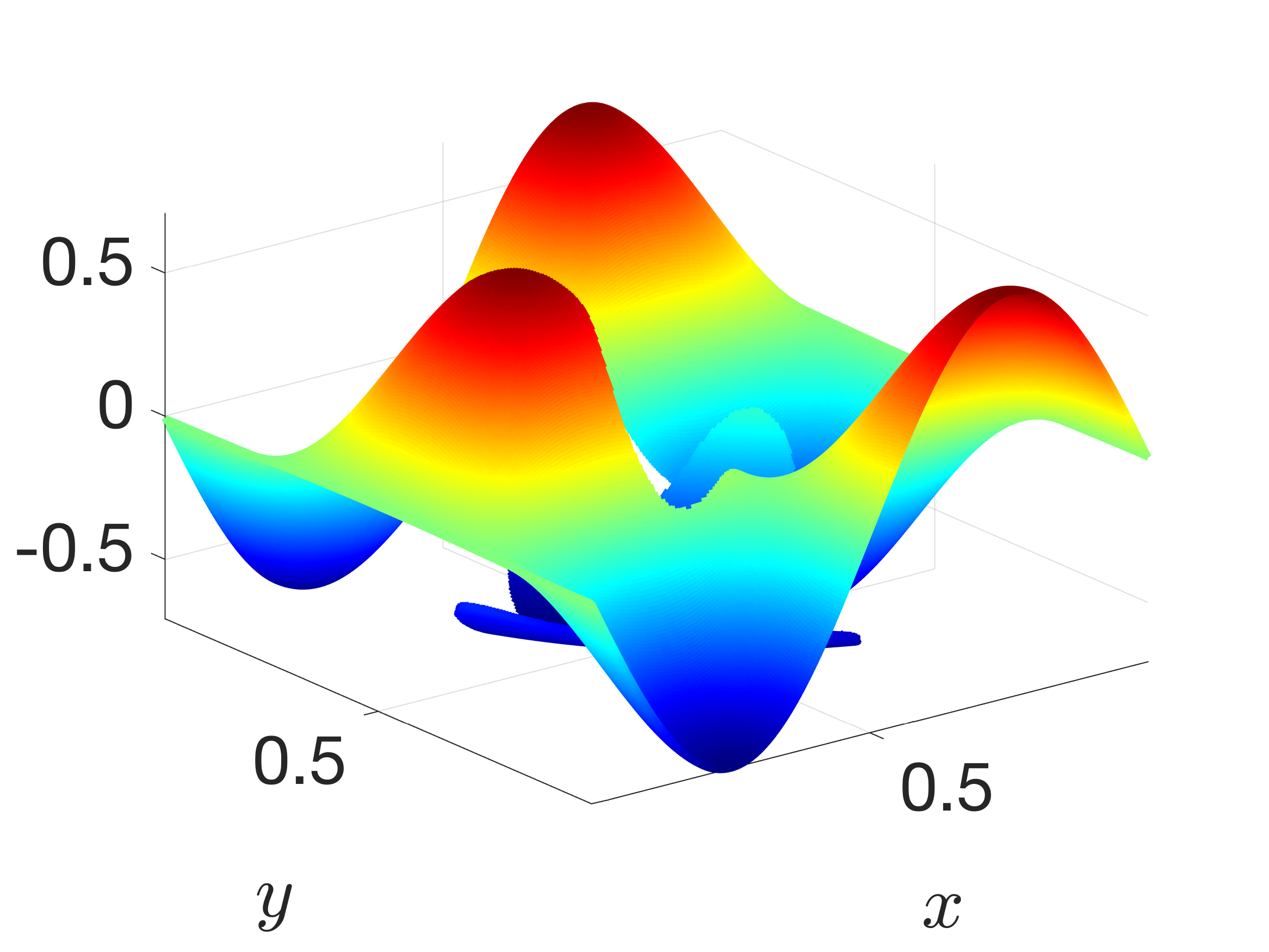}}{\footnotesize$E_z$}
	} 
	{\phantom{A}}
	\label{fig:starInterfacePEC}
	}
\subfigure[3-star]{
	\stackunder[5pt]{
	\stackunder[5pt]{\includegraphics[width=1.58in]{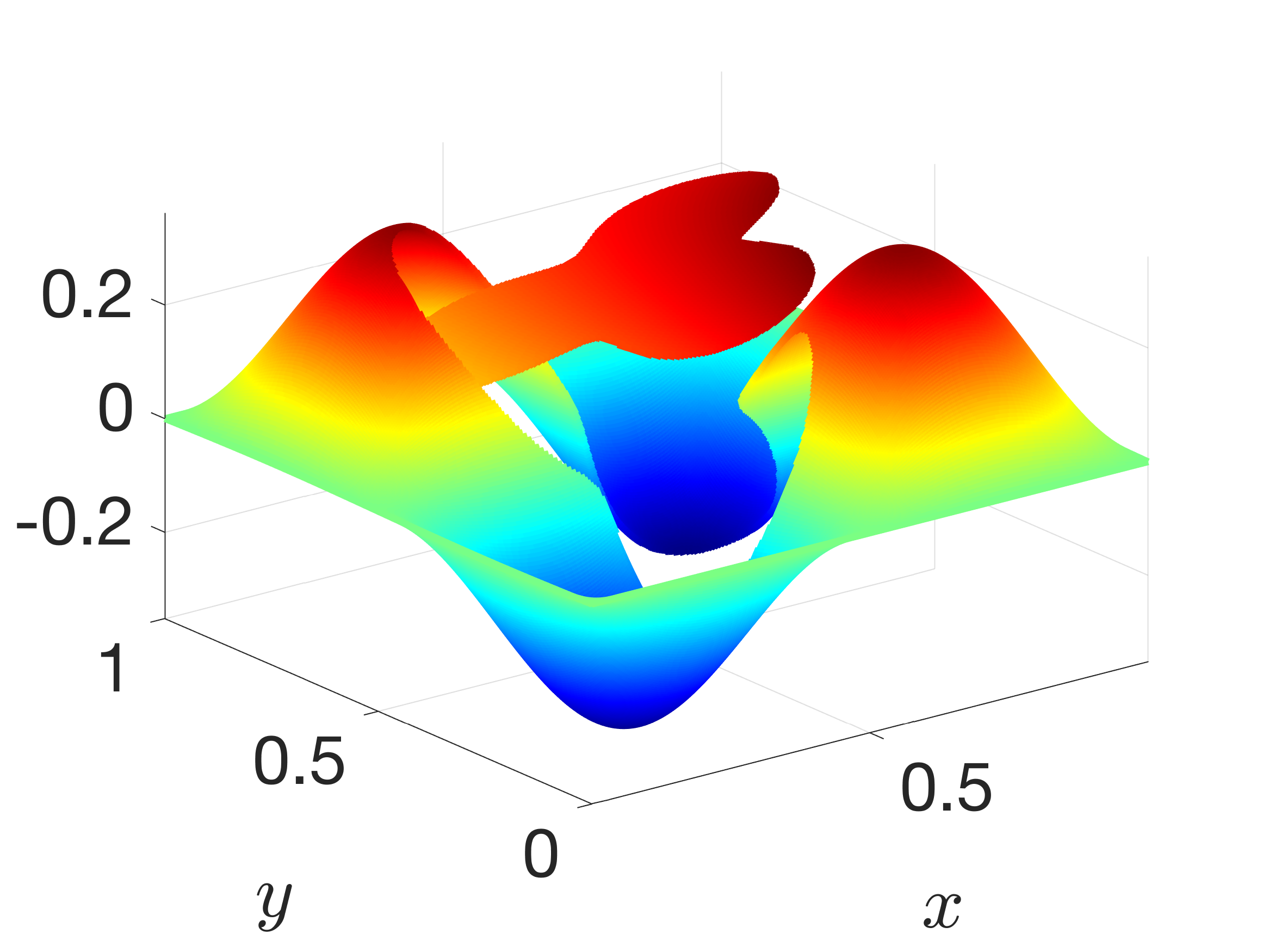}}{\footnotesize$H_x$}
	\stackunder[5pt]{	\includegraphics[width=1.58in]{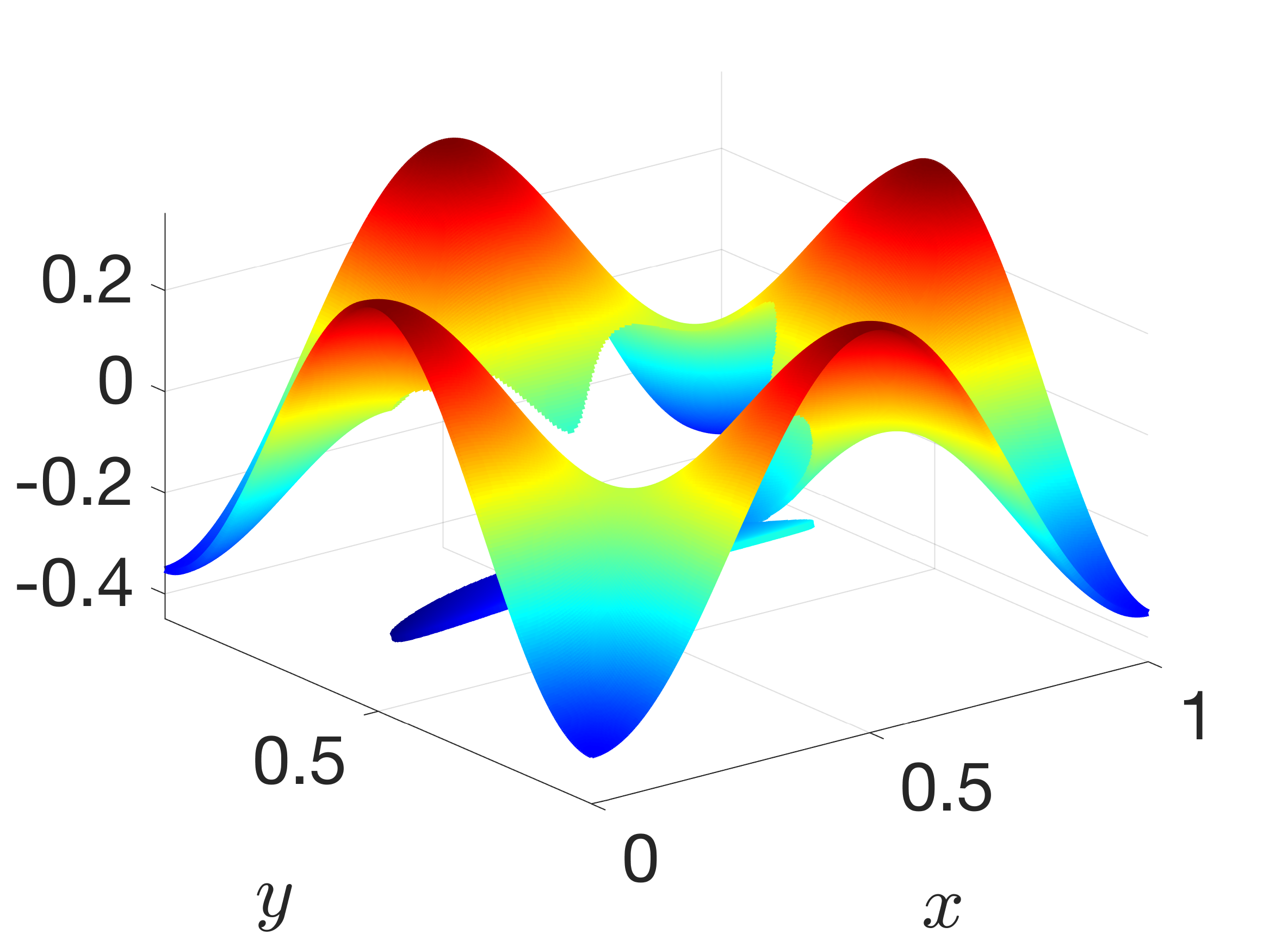}}{\footnotesize$H_y$}
	\stackunder[5pt]{	\includegraphics[width=1.58in]{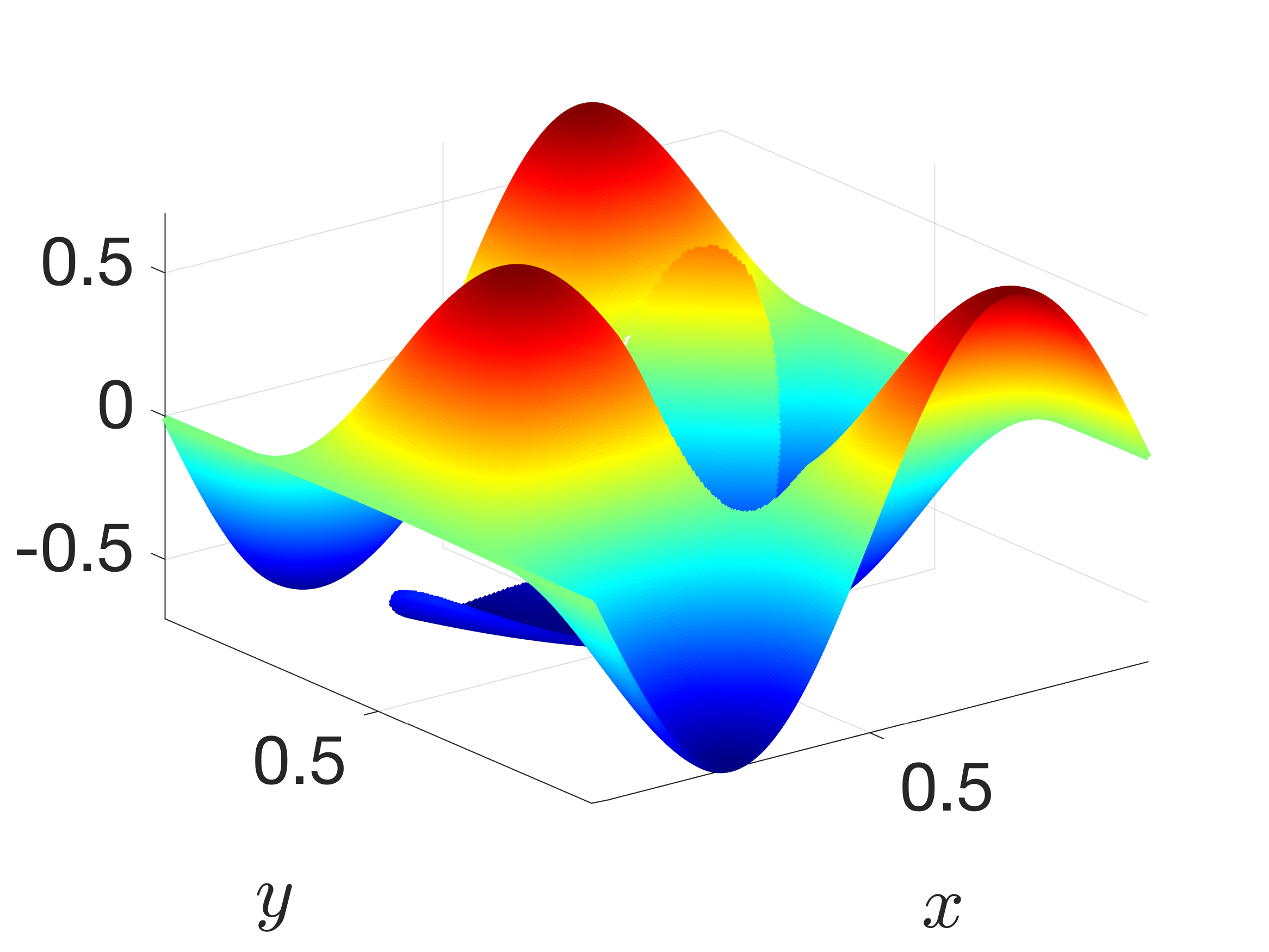}}{\footnotesize$E_z$}
	} 
	{\phantom{A}}
	\label{fig:triStarInterfacePEC}
	}
       \caption{The components $H_x$, $H_y$ and $E_z$ with $h = \tfrac{1}{336}$ for problems with a manufactured solution using the CFM-Yee scheme. The computed electric field and 
       magnetic field are shown respectively at $t=0.625$ and $t-\frac{\Delta t}{2}$.}
       \label{fig:plotInterfaceFields}
\end{figure}

\subsection{Stability Investigation : Long-Time Simulations} \label{sec:longTimeSimulations}
As mentioned in Remark~\ref{rem:stability}, 
	a rigorous stability analysis of CFM-FDTD schemes is out of reach for the moment. 
In this short subsection, 
	we therefore provide some numerical evidences on the stability of CFM-FDTD schemes for 
	a sufficiently small value of the penalization coefficient $c_f$.
We consider scattering of a dielectric cylinder problems, 
	and a problem with a manufactured solution and a 3-star interface. 
We use the CFM-Yee and the CFM-4$^{th}$ scheme. 
For both CFM-FDTD schemes, 
	the parameters remain the same as previously described. 
However, 
	we consider a larger time interval, 
	given by $I = [0, 25]$.
 \begin{figure} 
 \centering
  \subfigure[CFM-Yee]{ 
		\setlength\figureheight{0.33\linewidth} 
		\setlength\figurewidth{0.33\linewidth} 
		\tikzset{external/export next=false}
		\input{testNonMagneticLongTimeYee.tikz}
		}		
  \subfigure[CFM-4$^{th}$]{ 
		\setlength\figureheight{0.33\linewidth} 
		\setlength\figurewidth{0.33\linewidth} 
		\tikzset{external/export next=false}
		\input{testNonMagneticLongTime4th.tikz}
  }
  \caption{Evolution of the error in $L^2$-norm of $\mathbold{U} = [H_x,H_y,E_z]^T$ for a scattering of a dielectric cylinder problem with $\mu^+=\mu^-=1$ using the proposed CFM-FDTD schemes. The mesh grid size $\tfrac{1}{20}$, $\tfrac{1}{40}$ and $\tfrac{1}{80}$ correspond to respectively the black line, dotted blue line and dash-dotted magenta line.}
   \label{fig:longTimeSimulationNonMagnetic}
\end{figure}
%%%%%%%%%%%%%%%
 \begin{figure} 
 \centering
  \subfigure[CFM-Yee]{ 
		\setlength\figureheight{0.33\linewidth} 
		\setlength\figurewidth{0.33\linewidth} 
		\tikzset{external/export next=false}
		\input{testMagneticLongTimeYee.tikz}
		}		
  \subfigure[CFM-4$^{th}$]{ 
		\setlength\figureheight{0.33\linewidth} 
		\setlength\figurewidth{0.33\linewidth} 
		\tikzset{external/export next=false}
		\input{testMagneticLongTime4th.tikz}
  }
  \caption{Evolution of the error in $L^2$-norm of $\mathbold{U} = [H_x,H_y,E_z]^T$ for a scattering of a dielectric cylinder problem with $\mu^+=1$ and $\mu^-=2$ using the proposed CFM-FDTD schemes. The mesh grid size $\tfrac{1}{20}$, $\tfrac{1}{40}$ and $\tfrac{1}{80}$ correspond to respectively the black line, dotted blue line and dash-dotted magenta line.}
   \label{fig:longTimeSimulationMagnetic}
\end{figure}
%%%%%%%%%%%%%%%
 \begin{figure}
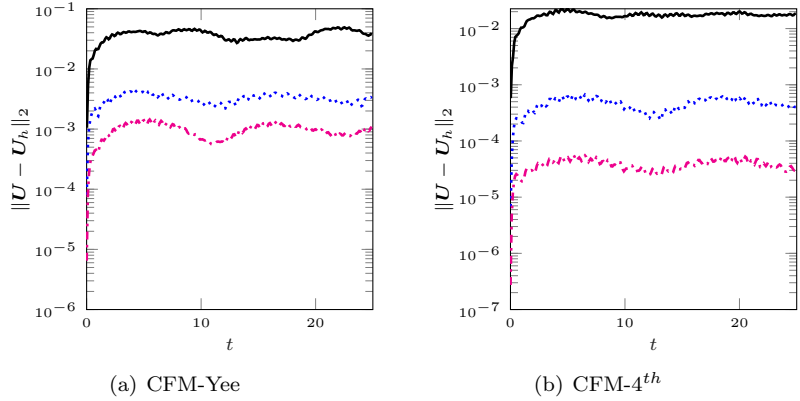
 
 \centering
  \subfigure[CFM-Yee]{ 
		\setlength\figureheight{0.33\linewidth} 
		\setlength\figurewidth{0.33\linewidth} 
		\tikzset{external/export next=false}
		\input{test3StarLongTimeYee.tikz}
		}		
  \subfigure[CFM-4$^{th}$]{ 
		\setlength\figureheight{0.33\linewidth} 
		\setlength\figurewidth{0.33\linewidth} 
		\tikzset{external/export next=false}
		\input{test3StarLongTime4th.tikz}
  }
  \caption{Evolution of the error in $L^2$-norm of $\mathbold{U} = [H_x,H_y,E_z]^T$ for a problem with a manufactured solution using a 3-star interface and the proposed CFM-FDTD schemes. The mesh grid size $\tfrac{1}{20}$, $\tfrac{1}{40}$ and $\tfrac{1}{80}$ correspond to respectively the black line, dotted blue line and dash-dotted magenta line.}
   \label{fig:longTimeSimulationManuSol}
\end{figure}
%%%%%%%%%%%%
Figure~\ref{fig:longTimeSimulationNonMagnetic}, 
	Figure~\ref{fig:longTimeSimulationMagnetic} 
	and Figure~\ref{fig:longTimeSimulationManuSol} 
	illustrate the evolution of the error in $L^2$-norm of $\mathbold{U} = [H_x,H_y,E_z]^T$ 
	for respectively a non-magnetic dielectric cylinder problem,
	a magnetic dielectric cylinder problem and a problem with a manufactured solution.
In all cases, 
	numerical results suggest that the proposed CFM-FDTD schemes are stable.

\section{Conclusions} \label{sec:conclusion}
In this work,
	we presented high-order FDTD schemes based on the Correction Function Method. 
The system of PDEs needed for the CFM
%	that models corrections 
%	that are needed to retain the order of a 
%	finite difference approximation close to the interface 
	was derived using Maxwell's equations with interface conditions.
The minimization problem based on a functional that is a square measure of the error associated with the 
	correction function's system of PDEs was also presented and solved.
%The functional that is a square measure of the error associated with the correction function's system of PDEs 
%	and the minimization problem are introduced.
Numerical examples showed that numerical solutions coming from CFM-FDTD schemes 
	were captured without spurious oscillation while exhibiting high-order convergence.
Moreover, 
	the accuracy of correction functions has been verified using high-order explicit jump conditions.
This showed that high-order jump conditions are implicitly enforced in the functional to minimize and therefore 
	need not be provided explicitly. 
Problems with a manufactured solution have shown that the proposed numerical strategy can 
	handle various geometries of the interface without significantly increasing the complexity 
	of the method.
Despite a lack of a rigorous stability analysis,
	long-time simulations have been performed and provided numerical evidences of 
	the stability of CFM-FDTD schemes.		
Future work will focus on the theoretical aspect of the CFM as well as an extension of 
	this strategy to 3-D problems. 

\section*{Acknowledgments}
The authors are grateful to Alexis Montoison of Polytechnique Montr\'{e}al for his help on Julia programming 
	language \cite{Bezanson2017}.
The authors also thank Dr. Jessica Lin and Dr. Gantumur Tsogtgerel of McGill University for their support.
The research of JCN was partially supported by the NSERC Discovery Program.
This preprint has not undergone peer review (when applicable) or any post-submission improvements or corrections. The Version of Record of this article is published in Journal of Scientific Computing, and is available online at https://doi.org/10.1007/s10915-022-01797-9.

%\appendix 
%\FloatBarrier
%\section{Convergence Tables} \label{sec:convTables}
%%***************************************
%\begin{table}%[h!]
%%\begin{center}
%\begin{adjustbox}{max width=1.5\textwidth,center}
%\renewcommand{\arraystretch}{1.5}
%%\footnotesize
%\scriptsize
%%\tiny
%\begin{tabular}{ccccc}
%\hline 
%h &    \multicolumn{2}{c}{CFM-Yee}&  \multicolumn{2}{c}{CFM-$4^{th}$}\\\cmidrule(lr){2-3}\cmidrule(lr){4-5}
% &   error  & order  & error & order   \\ \hline
%$\tfrac{1}{20}$  &  & - &   & -    \\ 
%$\tfrac{1}{28}$  & &  &  &  \\ 
%$\tfrac{1}{40}$ & &  &   &   \\ 
%$\tfrac{1}{52}$ & &  &   &     \\ 
%$\tfrac{1}{72}$  & &  &   &   \\  
%$\tfrac{1}{96}$  & &  &   &    \\   
%$\tfrac{1}{132}$  & &  &  &   \\  
%$\tfrac{1}{180}$  & &  &  &   \\ 
%$\tfrac{1}{244}$ & &  &   &   \\  
%$\tfrac{1}{336}$ & &     &  &       \\  \noalign{\smallskip} \hline
%\end{tabular}
%%\normalsize
%%\end{center}
%\end{adjustbox}
%\caption{Errors of $\mathbold{U}$ in $L^2$-norm and estimated orders for the scattering of a dielectric cylinder problem with $\mu^-=1$, $\epsilon^-= 2.25$, $\mu^+ = 1$ and $\epsilon^+ = 1$ using the proposed CFM-FDTD schemes.}
%\label{tab:scatteringCylinderNonMagneticConvTab}
%\end{table}
%***************************************

\section*{References}

\bibliography{references}

\end{document}